\documentclass[fleqn]{article}
\usepackage{makeidx}
\usepackage{amsmath}
\usepackage{amssymb}
\usepackage{indentfirst}
\usepackage{amsfonts}

\usepackage{amsthm}
\newtheorem{theorem}{Theorem}
\newtheorem*{acknowledgement}{Acknowledgement}

\newtheorem{corollary}[theorem]{Corollary}

\newtheorem{definition}[theorem]{Definition}
\newtheorem{example}[theorem]{Example}

\newtheorem{lemma}[theorem]{Lemma}

\newtheorem{proposition}[theorem]{Proposition}
\newtheorem{remark}[theorem]{Remark}

\input{tcilatex}
\roman{enumii}
\numberwithin{equation}{section}

\begin{document}

\title{Galois extensions for coquasi-Hopf algebras}
\author{ADRIANA BALAN \\
Faculty of Applied Sciences, University Politehnica of Bucharest, \\
313 Splaiul Independen\c{t}ei, 060042 Bucharest, Romania\\
e-mail: asteleanu@yahoo.com}
\date{April, 2008}
\maketitle

\begin{abstract}
The notions of Galois and cleft extensions are generalized for coquasi-Hopf
algebras. It is shown that such an extension over a coquasi-Hopf algebra is
cleft if and only if it is Galois and has the normal basis property. A
Schneider type theorem (\cite{Schneider90a}) is proven for coquasi-Hopf
algebras with bijective antipode. As an application, we generalize
Schauenburg's bialgebroid construction for coquasi-Hopf algebras.
\end{abstract}

\section{Introduction}

The study of Hopf Galois extensions started last century with the papers of
Chase, Harrison and Rosenberg (\cite{Chase65})\ and of Chase and Sweedler (%
\cite{Chase69}). Later, their definition was improved by Kreimer and
Takeuchi (\cite{Kreimer81}) and knew since a continuous development, mainly
because their relation to different areas of mathematics. But in the last
decade, examples of extensions which were not Galois but behaved like such
had appeared. An explanation was necessary, and it became soon clear that
this was possible only by generalization. The replacement of the Hopf
algebra by a coalgebra (or more generally by a coring) has led to the notion
of a Galois extension by a coalgebra, first formulated by Brz\`{e}zinski and
Hajac (\cite{Brzezinski99}). Another generalization was obtained considering
Galois extensions over a coacting bialgebroid (over a non-commutative ground
ring) (\cite{Kadison03}, \cite{Bohm04}).

All structures cited above are generalizations of bialgebras or Hopf
algebras. Another such objects are the (co)quasi-Hopf algebras. They have
been introduced by Drinfeld (\cite{Drinfeld90}), respectively by Majid (\cite%
{Majid92}) and have lately attracted much attention in both mathematics and
physics (\cite{Altschuler92}, \cite{Majid95}). So it is natural to see if it
is possible to generalize the Galois theory also to the case of coquasi-Hopf
algebras.

The definition of a coquasi-Hopf algebra $H$ ensures that the category of
right $H$-comodules $\mathcal{M}^{H}$ is monoidal, with usual tensor product
over the base field. The difference between a coquasi-Hopf algebra and a
Hopf algebra is that the associativity of tensor product in the monoidal
category does not coincide with the usual associativity of tensor product in
the category of vector spaces. Consequently, the multiplication of a
coquasi-Hopf algebra is no longer associative, but associative up to
conjugation by an invertible element $\omega \in (H\otimes H\otimes H)^{\ast
}$ (the reassociator). But is this main feature of coquasi-Hopf algebras,
namely the monoidallity of corepresentations, which made possible
generalizations of major properties from Hopf algebras (the existence and
uniqueness of integrals, the Nichols-Zoeller Theorem, construction of the
Drinfeld double, etc.). Hence it seems natural to continue with the Galois
theory for coquasi-Hopf algebras.

The present paper begins with a short review of the known results about
coquasi-Hopf algebras, their categories of comodules and about algebras and
modules within these monoidal categories mentioned above. As for Hopf
algebras, to each right comodule algebra $A$ (which is an algebra in the
monoidal category of right comodules) one can assign a pair of adjoint
functors, namely the functor of coinvariants and the induced functor. Our
purpose is to generalize their well-known properties from Hopf algebras to
coquasi-Hopf algebras.

In the second part it is defined the notion of Galois extension. A right
comodule algebra $A$ is a Galois extension over its coinvariants ring $%
A^{coH}$(which is associative, although $A$ fails to be) if a certain map is
bijective. This is a natural generalization of the author's previous paper (%
\cite{Balan06}), where only finite dimensional quasi-Hopf algebras were
considered. It should be noticed that this definition for Galois extensions
works only for coquasi-Hopf algebras, as it involves the presence of the
antipode. Although this may look restrictive, we shall see that this
definition for the Galois map allows us to recover all principal results
from the classical Hopf-Galois theory. A Galois extension is invariant to
any gauge transformation. As an example of Galois extension, we take a group
algebra, and view it as a coquasi-Hopf algebra by a 3-cocycle. Then a
comodule algebra is a graded space with a multiplication non-necessarily
associative, which is Galois over its invariants if and only if it is
strongly graded. This was known from long time for Hopf algebras (\cite%
{Ulbrich81}). Moreover, any crossed product (\cite{Balan07a}) coming from an
associative algebra endowed with a 2-cocycle and a weak action is also a
Galois extension.

Recall that in the Hopf algebra case, the functor of coinvariants is a $Hom$
functor, and the Galois map is just the evaluation for a certain relative
Hopf module. We show that these results, slightly modified by the presence
of a twist, hold also in the coquasi-case. We give thus an explanation for
the formula of the Galois map from Definition \ref{def Galois}. Although the
results are the same, it is much more difficult to obtain them. The
structure of the relative Hopf module $A\otimes H$ (which is the link
between the Galois map $can$ and the adjunction of categories $\mathcal{M}%
_{B}\overset{(-)\otimes _{B}A}{\underset{(-)^{coH}}{\mathcal{%
\rightleftarrows }}}\mathcal{M}_{A}^{H}$) is not obvious. The classical
formulas do not work anymore, and an isomorphism is required in order to get
the desired structure by transport.

We introduce next the notion of a cleft extension. As this involves the
convolution product (which is no longer associative), the invertibility of
the cleaving map has to be translated now in relations (\ref%
{convolutiedeltagama}), (\ref{convolutiegamabetadelta}) involving the
antipode and the linear maps $\alpha $, $\beta $.

As a generalization of theorems of Doi and Takeuchi (\cite{Doi86}), and
Blattner and Montgomery (\cite{Blattner89}), we obtain the first main result
of this paper, namely the equivalence between cleft extensions, and Galois
extensions with the normal basis property.

The second main part of this section concerns the equivalence between the
category of relative Hopf modules and modules over the subalgebra of
coinvariants. It starts with an analogue of the Schneider%
\'{}%
s imprimitivity theorem of \cite{Schneider90a}. A key problem in the proof
is how to show that the bijectivity of the Galois map implies the
bijectivity of the corresponding map for any relative right $(A,H)$-Hopf
module. As $A$ is not an associative algebra, this is not obvious and
requires some special considerations about the tensor product over the
algebra $A$ in the monoidal category of right comodules (Lemma \ref{canM
bijectiv}). The proof of the theorem uses the Five Lemma applied twice to
some commutative diagrams, but unlike the Hopf algebra case, the
commutativity of those is not an easy fact and requires special attention
care when dealing with the reassociator $\omega $ and of the elements $%
\alpha $ and $\beta $ (from the definition of the antipode).

Next, we prove a coquasi-version of the affineness criteria for affine
algebraic groups schemes, where the surjectivity of the Galois map of the
extension is related to relative injectivity of the $H$-comodule $A$ and to
the equivalence between the category of relative Hopf modules and modules
over the subalgebra of coinvariants.

In the last section, we generalize Schauenburg's bialgebroid construction.
This is an illustration of how the Galois theory, combined with monoidally
arguments can raise to new structures.

\section{Preliminaries}

In this section we recall some definitions and results and fix notations.
Throughout the paper we work over some base field $k$. Tensor products,
algebras, linear spaces, etc. will be over $k$. Unadorned $\otimes $ means $%
\otimes _{k}$. We shall use dots to indicate the module or comodule
structure on the tensor product. An introduction to the study of
quasi-bialgebras and quasi-Hopf algebras and their duals
(coquasi-bialgebras, respectively coquasi-Hopf algebras) can be found in 
\cite{Majid95}. A good reference for monoidal categories is \cite{Kassel95},
while actions of monoidal categories are exposed in \cite{Pareigis77}, \cite%
{Pareigis77II}.

\begin{definition}
A \textbf{coquasi-bialgebra} $(H,m,u,\omega ,\Delta ,\varepsilon )$ is a
coassociative coalgebra $(H,\Delta ,\varepsilon )$ together with coalgebra
morphisms: the multiplication $m:H\otimes H\longrightarrow H$ (denoted $%
m(h\otimes g)=hg$), the unit $u:\Bbbk \longrightarrow H$ (denoted $%
u(1)=1_{H} $), and a convolution invertible element $\omega \in (H\otimes
H\otimes H)^{\ast }$ such that: 
\begin{eqnarray}
h_{1}(g_{1}k_{1})\omega (h_{2},g_{2},k_{2}) &=&\omega
(h_{1},g_{1},k_{1})(h_{2}g_{2})k_{2}  \label{asociat multipl} \\
1_{H}h &=&h1_{H}=h \\
\omega (h_{1},g_{1},k_{1}l_{1})\omega (h_{2}g_{2},k_{2},l_{2}) &=&\omega
(g_{1},k_{1},l_{1})\omega (h_{1},g_{2}k_{2},l_{2})\omega (h_{2},g_{3},k_{3})
\label{cocycle omega} \\
\omega (h,1_{H},g) &=&\varepsilon (h)\varepsilon (g)
\end{eqnarray}%
hold for all $h,g,k,l\in H$.
\end{definition}

As a consequence, we have also $\omega (1_{H},h,g)=\omega
(h,g,1_{H})=\varepsilon (h)\varepsilon (g)$ for each $g,h\in H$.

\begin{definition}
A \textbf{coquasi-Hopf algebra} is a coquasi-bialgebra $H$ endowed with a
coalgebra antihomomorphism $S:H\longrightarrow H$ (the antipode) and with
elements $\alpha $, $\beta \in H^{\ast }$ satisfying 
\begin{eqnarray}
S(h_{1})\alpha (h_{2})h_{3} &=&\alpha (h)1_{H}  \label{SalfaId} \\
h_{1}\beta (h_{2})S(h_{3}) &=&\beta (h)1_{H}  \label{IdbetaS} \\
\omega (h_{1}\beta (h_{2}),S(h_{3}),\alpha (h_{4})h_{5}) &=&\omega
^{-1}(S(h_{1}),\alpha (h_{2})h_{3}\beta (h_{4}),S(h_{5}))=\varepsilon (h)
\label{omega anihileaza S}
\end{eqnarray}%
for all $h\in H$.
\end{definition}

These relations imply also $S(1_{H})=1_{H}$ and $\alpha (1_{H})\beta
(1_{H})=1$, so by rescaling $\alpha $ and $\beta $, we may assume that $%
\alpha (1_{H})=1$ and $\beta (1_{H})=1$. The antipode is unique up to a
convolution invertible element $U\in H^{\ast }$: if $(S^{\prime },\alpha
^{\prime },\beta ^{\prime })$ is another triple with the above properties,
then according to \cite{Majid95} we have 
\begin{equation}
S^{\prime }(h)=U(h_{1})S(h_{2})U^{-1}(h_{3}),\qquad \alpha ^{\prime
}(h)=U(h_{1})\alpha (h_{2}),\qquad \beta ^{\prime }(h)=\beta
(h_{1})U^{-1}(h_{2})  \label{change antipode coquasi}
\end{equation}%
for all $h\in H$.

We shall use in this paper the monoidal structure of the right $H$-comodule
category $\mathcal{M}^{H}$ and of the left $H$-comodule category$\ {}^{H}%
\mathcal{M}$: the tensor product is over the base field and the comodule
structure (left or right) of the tensor product is the codiagonal one. The
reassociators are 
\begin{eqnarray*}
\phi _{U,V,W} &:&(U\otimes V)\otimes W\longrightarrow U\otimes (V\otimes W)
\\
\phi _{U,V,W}((u\otimes v)\otimes w) &=&u_{0}\otimes (v_{0}\otimes
w_{0})\omega (u_{1},v_{1},w_{1})
\end{eqnarray*}%
for $u\in U$, $v\in V$, $w\in W$ and $U,V,W\in \mathcal{M}^{H}$,
respectively 
\begin{eqnarray*}
\phi _{U,V,W} &:&(U\otimes V)\otimes W\longrightarrow U\otimes (V\otimes W)
\\
\phi _{U,V,W}((u\otimes v)\otimes w) &=&\omega
^{-1}(u_{-1},v_{-1},w_{-1})u_{0}\otimes (v_{0}\otimes w_{0})
\end{eqnarray*}%
for $u\in U$, $v\in V$, $w\in W$ and $U,V,W\in {}^{H}\mathcal{M}$.

Together with a coquasi-Hopf algebra with bijective antipode $H=(H,\Delta
,\varepsilon ,m,1_{H},\omega ,S,\alpha ,\beta )$, we also have $H^{op}$, $%
H^{cop}$, and $H^{op,cop}$ as coquasi-Hopf algebras, where "op" means
opposite multiplication and "cop" means opposite comultiplication. The
coquasi-Hopf structures are obtained by putting $\omega _{cop}=\omega ^{-1}$%
, $\omega _{op}=(\omega ^{-1})^{321}$, $\omega _{op,cop}=\omega ^{321}$, $%
S_{op}=S_{cop}=(S_{op,cop})^{-1}=S^{-1}$, $\alpha _{cop}=\beta S^{-1}$,$\,\
\alpha _{op}=\alpha S^{-1}$, $\alpha _{op,cop}=\beta $, $\beta _{cop}=\alpha
S^{-1}$, $\beta _{op}=\beta S^{-1}$ and $\beta _{op,cop}=\alpha $. Here $%
\omega ^{321}(h,g,k)=\omega (k,g,h)$.

For $H$ a coquasi-bialgebra, the linear dual $H^{\ast }=Hom(H,\Bbbk )$
becomes an associative algebra with multiplication given by the convolution
product%
\begin{equation}
(h^{\ast }g^{\ast })(h)=h^{\ast }(h_{1})g^{\ast }(h_{2})\qquad \forall h\in
H\quad \text{\c{s}i}\quad h^{\ast },g^{\ast }\in H^{\ast }
\label{convolutie pe hrond}
\end{equation}%
and unit $\varepsilon $. This algebra is acting on $H$ by the formulas:%
\begin{equation}
h^{\ast }\rightharpoonup h=h_{1}h^{\ast }(h_{2})\text{,}\qquad
h\leftharpoonup h^{\ast }=h^{\ast }(h_{1})h_{2}
\label{actiunea slaba a lui Hrond pe Hrond*}
\end{equation}%
for any $h^{\ast }\in H^{\ast }$, $h\in H$.

Now, recall from \cite{Panaite97Stefan} the following: for $\tau \in
(H\otimes H)^{\ast }$ a convolution invertible map such that $\tau
(1,h)=\tau (h,1)=\varepsilon (h)$ for all $h\in H$ ($\tau $ is called a
twist or a gauge transformation), one can define a new structure of
coquasi-Hopf algebra on $H$, denoted $H_{\tau }$, by taking 
\begin{eqnarray}
h\cdot _{\tau }g &=&\tau (h_{1},g_{1})h_{2}g_{2}\tau ^{-1}(h_{3},g_{3})
\label{multiplic Htau} \\
\omega _{\tau }(h,g,k) &=&\tau (g_{1},k_{1})\tau (h_{1},g_{2}k_{2})\omega
(h_{2},g_{3},k_{3})\tau ^{-1}(h_{3}g_{4},k_{4})\tau ^{-1}(h_{4},g_{5})
\label{asociator Htau} \\
\alpha _{\tau }(h) &=&\tau ^{-1}(S(h_{1}),\alpha (h_{2})h_{3})
\label{alfa Htau} \\
\beta _{\tau }(h) &=&\tau (h_{1}\beta (h_{2}),S(h_{3}))  \label{beta Htau}
\end{eqnarray}%
for all $h,g,k\in H$, and keeping the unit, the comultiplication, the counit
and the antipode unchanged.

\begin{remark}
\label{monoidal isomorfism gauge}There is a monoidal isomorphism $\mathcal{M}%
^{H}\simeq \mathcal{M}^{~H_{\tau }\text{, }}$which is the identity on
objects and on morphisms, with monoidal structure given by $V\otimes
W\longrightarrow V\otimes W$, $v\otimes w\longrightarrow v_{0}\otimes
w_{0}\tau ^{-1}(v_{1},w_{1})$, where $v\in V$, $w\in W$ and $V$, $W\in 
\mathcal{M}^{H}$.
\end{remark}

In \cite{Bulacu99}, it was constructed a twist \textbf{$f$}$\in (H\otimes
H)^{\ast }$ which controls how far is $S$ from a anti-algebra morphism:%
\begin{equation}
\mathbf{f}(h_{1},g_{1})S(h_{2}g_{2})=S(g_{1})S(h_{1})\mathbf{f}(h_{2},g_{2})
\label{twist f}
\end{equation}%
If we denote%
\begin{eqnarray}
p(h,g) &=&\omega (S(g_{2}),S(h_{2}),h_{4})\omega
^{-1}(S(g_{1})S(h_{1}),h_{5},g_{4})\alpha (h_{3})\alpha (g_{3})
\label{formula p} \\
q(h,g) &=&\omega (h_{1}g_{1},S(g_{5}),S(h_{4}))\omega
^{-1}(h_{2},g_{2},S(g_{4}))\beta (h_{3})\beta (g_{3})  \label{formula q}
\end{eqnarray}%
then the twist $\mathbf{f}$ is given by%
\begin{equation*}
\mathbf{f}(h,g)=\omega
^{-1}(S(g_{1})S(h_{1}),h_{3}g_{3},S(h_{5}g_{5}))p(h_{2},g_{2})\beta
(h_{4}g_{4})
\end{equation*}%
We have also that%
\begin{eqnarray}
\mathbf{f}(h_{1},g_{1})\alpha (h_{2}g_{2}) &=&p(h,g)
\label{f*alfa=gama, beta*f-1=delta} \\
\beta (h_{1}g_{1})\mathbf{f}^{(-1)}(h_{2},g_{2}) &=&q(h,g)
\label{beta*f-1=delta} \\
p(h_{1},S(h_{3}))\beta (h_{2}) &=&\alpha S(h)  \label{relatie p} \\
\mathbf{f}(h_{1},S(h_{3}))\beta (h_{2}) &=&\alpha S(h)  \label{relatie f} \\
\mathbf{f}^{(-1)}(S^{-1}(g_{1}),S^{-1}(h_{1}))\omega ^{-1}(g_{4},\alpha
S^{-1}(g_{3})S^{-1}(g_{2}),S^{-1}(h_{2})) &=&\mathbf{f}%
(g_{5},S^{-1}(h_{1}g_{1}))  \label{relatie UL pR h} \\
&&\omega ^{.1}(h_{2},g_{2}\beta (g_{3}),S(g_{4}))  \notag
\end{eqnarray}%
where in the last formula we assumed the bijectivity of the antipode.
Relations (\ref{f*alfa=gama, beta*f-1=delta}) and (\ref{beta*f-1=delta}) are
from \cite{Bulacu02co}, (\ref{relatie p}) is an easy consequence of the
formula of $p$, (\ref{relatie f}) follows immediately from (\ref{relatie p}%
), while for (\ref{relatie UL pR h}) we use (\ref{asociator Htau})\ and the
fact that the associator $\omega _{\mathbf{f}}$ for the twisted
coquasi-bialgebra $H_{\mathbf{f}}$ is $\omega _{\mathbf{f}}(h,g,k)=\omega
(S(k),S(g),S(h))$, $\forall $ $h,g,k\in H$. If the antipode is bijective,
then by passing from $H$ to $H^{op}$ we obtain a new twist $\widetilde{%
\mathbf{f}}\in (H\otimes H)^{\ast }$ \cite{BulacuChirita01}, given by%
\begin{equation}
\widetilde{\mathbf{f}}(h,g)=\mathbf{f}(S^{-1}(g),S^{-1}(h))  \label{twist h}
\end{equation}%
which satisfies 
\begin{equation}
\widetilde{\mathbf{f}}%
(h_{1},g_{1})S^{-1}(h_{2}g_{2})=S^{-1}(g_{1})S^{-1}(h_{1})\widetilde{\mathbf{%
f}}(h_{2},g_{2})  \label{hdeltaS-1(a)h-1=S-1tensorS-1(delta cop(a))}
\end{equation}%
for any $h,g\in H$. The corresponding reassociator will be $\omega _{%
\widetilde{\mathbf{f}}}(h,g,k)=\omega (S^{-1}(k),S^{-1}(g),S^{-1}(h))$. This
twist will appear later.

\begin{definition}
(\cite{Bulacu00}) A right $H$-comodule algebra $A$ is an algebra in the
monoidal category $\mathcal{M}^{H}$. This means $(A,\rho _{A})$ is a right $%
H $-comodule, we have a multiplication map $m_{A}:A\otimes A\longrightarrow
A $, denoted $m_{A}(a\otimes b)=ab$, for $a,b\in A$, and a unit map $%
u_{A}:\Bbbk \longrightarrow A$, where we put $u_{A}(1)=1_{A}$, which are
both $H$-colinear, such that 
\begin{equation}
(ab)c=a_{0}(b_{0}c_{0})\omega (a_{1},b_{1},c_{1})  \label{asoc comod alg}
\end{equation}%
holds for any $a,b,c\in A$.
\end{definition}

Similarly we may define a left $H$-comodule algebra as an algebra in $^{H}%
\mathcal{M}$. Notice that $A$ is a right $H$-comodule algebra if and only if 
$A^{op}$ is a left $H^{op,cop}$-comodule algebra.

\begin{definition}
(\cite{Bulacu00}) For $A$ a right $H$-comodule algebra, we may define the
notion of right module over $A$ in the category $\mathcal{M}^{H}$.\
Explicitly, this is a right $H$-comodule $(M,\rho _{M})$, endowed with a
right $A$-action, denoted $\mu _{M}(m,a)=ma$, such that 
\begin{eqnarray*}
(ma)b &=&m_{0}(a_{0}b_{0})\omega (m_{1},a_{1},b_{1}) \\
m1_{A} &=&m \\
\rho _{M}(ma) &=&m_{0}a_{0}\otimes m_{1}a_{1}
\end{eqnarray*}%
hold for all $m\in M$, $a,b\in A$. The category of such objects, with
morphisms the right $H$-colinear maps which respect the $A$-action, is
called the \textbf{category of relative right }$(H,A)$\textbf{-Hopf modules}
and denoted $\mathcal{M}_{A}^{H}$.
\end{definition}

In the same way, we may define the category of left relative Hopf modules $%
_{A}\mathcal{M}^{H}$ for $A$ a right $H$-comodule algebra. If $A$ is a left $%
H$-comodule algebra we can define similarly the categories $_{A}^{H}\mathcal{%
M}$ and $^{H}\mathcal{M}_{A}$. For later use, remark that the following
categories are isomorphic:%
\begin{equation}
_{A}^{H}\mathcal{M}\simeq \mathcal{M}_{A^{op}}^{H^{op,cop}}
\label{left-right}
\end{equation}
for any $A$ a left $H$-comodule algebra (\cite{Bulacu00}).

\begin{remark}
\label{twist comodule algebra}It was proven in \cite{Bulacu00} that if $\tau 
$ is a twist on $H$, then the formula 
\begin{equation}
a\cdot _{\tau }b=a_{0}b_{0}\tau ^{-1}(a_{1},b_{1})  \label{multiplic Atau}
\end{equation}%
for all $a,b\in A$ defines a new multiplication such that $A$, with this new
multiplication (denoted $A_{\tau ^{-1}}$) becomes a right $H_{\tau }$%
-comodule algebra.\ It is easy to see that the isomorphism of Remark \ref%
{monoidal isomorfism gauge} sends the algebra $A$ of the monoidal category $%
\mathcal{M}^{H}$ exactly to the algebra $A_{\tau ^{-1}}$ in $\mathcal{M}%
^{H_{\tau }}$. But $\mathcal{M}^{H}$ and $\mathcal{M}^{H_{\tau }}$ are
monoidally isomorphic, therefore the categories of right relative Hopf
modules $\mathcal{M}_{A}^{H}$ and $\mathcal{M}_{A_{\tau ^{-1}}}^{H_{\tau }}$
will also be isomorphic.
\end{remark}

Let $A$ be a right $H$-comodule algebra. Consider the space of coinvariants%
\begin{equation*}
B=A^{coH}=\{a\in A\left\vert \rho _{A}(a)=a\otimes 1_{H}\right. \}
\end{equation*}%
It is immediate that this is an associative algebra with unit and
multiplication induced by the unit and the multiplication of $A$.

Now for each $M\in \mathcal{M}_{A}^{H}$, denote $M^{coH}=\{m\in M\mathcal{%
\mid }\rho _{M}(m)=m\otimes 1_{H}\}$. Then $M^{coH}$ becomes naturally a
right $B$-module, so we get the coinvariant functor 
\begin{equation*}
\mathcal{M}_{A}^{H}\overset{(-)^{coH}}{\mathcal{\longrightarrow }}\mathcal{M}%
_{B}
\end{equation*}%
Notice also the natural isomorphism 
\begin{equation}
Hom_{A}^{H}(A,M)\simeq M^{coH}  \label{coH=hom}
\end{equation}%
for any $M\in \mathcal{M}_{A}^{H}$. Conversely, for $N\in \mathcal{M}_{B}$,
we have $N\otimes _{B}A\in \mathcal{M}_{A}^{H}$ by 
\begin{eqnarray*}
\rho (n\otimes _{B}a) &=&n\otimes _{B}a_{0}\otimes a_{1} \\
(n\otimes _{B}a)b &=&n\otimes _{B}ab
\end{eqnarray*}%
As in the classical Hopf algebra case, we obtain the following:

\begin{proposition}
\label{functorul indus}The induced functor $(-)\otimes _{B}A$ is a left
adjoint for the functor of coinvariants $(-)^{coH}$:%
\begin{equation*}
\mathcal{M}_{B}\overset{(-)\otimes _{B}A}{\underset{(-)^{coH}}{\mathcal{%
\rightleftarrows }}}\mathcal{M}_{A}^{H}
\end{equation*}
\end{proposition}

\begin{proof}
Straightforward. For later use, we mention the adjunction morphisms:%
\begin{gather*}
\varepsilon _{M}:M^{coH}\otimes _{B}A\longrightarrow M,\qquad \varepsilon
_{M}(m\otimes _{B}a)=ma \\
u_{N}:N\longrightarrow (N\otimes _{B}A)^{coH},\qquad u_{N}(n)=n\otimes
_{B}1_{A}
\end{gather*}%
for each $N\in \mathcal{M}_{B}$ and $M\in \mathcal{M}_{A}^{H}$. Using the
isomorphism form relation \ref{coH=hom}, we get that the counit of the
adjunction is simply the evaluation.
\end{proof}

Similarly we could define the left version of the adjunction between the
induced and the coinvariant functor, namely $\mathcal{M}_{B}\overset{%
A\otimes _{B}(-)}{\underset{(-)^{coH}}{\mathcal{\rightleftarrows }}}{}_{A}%
\mathcal{M}^{H}$.

In the next section we shall see necessary and sufficient conditions for
these adjunctions to be equivalences.

\section{Galois extensions}

Let $H$ be a coquasi-Hopf algebra with antipode $S$ and $A$ a right $H$%
-comodule algebra. Denote as before $B=A^{coH}$.

\begin{definition}
\label{def Galois}The extension $B\subseteq A$ is $(H,S)$-\textbf{Galois} if
the map $can_{S}:A\otimes _{B}A\longrightarrow A\otimes H$, given by 
\begin{equation}
a\otimes _{B}b\longrightarrow a_{0}b_{0}\otimes b_{4}\omega
^{-1}(a_{1},b_{1}\beta (b_{2}),S(b_{3}))  \label{can}
\end{equation}%
is bijective.
\end{definition}

\begin{remark}
\rm%
(1) Although $A$ is not an associative algebra, we still keep the expression
"extension $B\subseteq A$".

(2) Recall that for coquasi-Hopf algebras the antipode is unique up to
conjugation to an invertible element. Therefore we need to check what is
happening if we change $S$. Consider another triple $(S^{\prime },\alpha
^{\prime },\beta ^{\prime })$ given by a convolution invertible element $%
U\in H^{\ast }$, as in (\ref{change antipode coquasi}). Then we have%
\begin{eqnarray*}
can_{S^{\prime }}(a\otimes _{B}b) &=&a_{0}b_{0}\otimes \omega
^{-1}(a_{1},b_{1}\beta ^{\prime }(b_{2}),S^{\prime }(b_{3}))b_{4} \\
&=&a_{0}b_{0}\otimes \omega ^{-1}(a_{1},b_{1}\beta
(b_{2})U^{-1}(b_{3}),U(b_{4})S(b_{5})U^{-1}(b_{6}))b_{7} \\
&=&a_{0}b_{0}\otimes \omega ^{-1}(a_{1},b_{1}\beta
(b_{2}),S(b_{3}))U^{-1}(b_{4})b_{5}
\end{eqnarray*}%
for every $a,b\in A$. If we define the linear map $\psi _{U}:A\otimes
H\longrightarrow A\otimes H$, $a\otimes h\longrightarrow a\otimes
U(h_{1})h_{2}$, it is easy to see that this is bijective with inverse $%
a\otimes h\longrightarrow a\otimes U^{-1}(h_{1})h_{2}$ and that $%
can_{S^{\prime }}=\psi _{U}\circ can_{S}$; therefore the two Galois maps
will be simultaneously bijective. In the sequel, we shall fix the antipode $%
S $ and the elements $\alpha $, $\beta $, such that $\alpha (1)=\beta (1)=1$%
, and write simply $can$.
\end{remark}

In case of a Hopf algebra, the coassociator $\omega $ and the linear map $%
\beta $ vanish, and we recover the usual definition of the Galois map. But
unlike the Hopf case, notice this time the presence of the antipode in the
formula of $can$, which implies that this definition is possible only for
coquasi-Hopf algebras, not also for coquasi-bialgebras. However, we shall
see that this definition for the Galois map allows us to recover all
principal results from the classical Hopf-Galois theory. In \cite{Masuoka03}%
, Masuoka uses the classical definition of the $can$ map, $a\otimes
_{B}b\longrightarrow ab_{0}\otimes b_{1}$, to show that a certain extension
is Galois over a given coquasi-Hopf algebra (which is a bicrossed product
associated to some cocyle data). It is only a matter of computation to see
that in the quoted case, the formula (\ref{can}) reduces to $a\otimes
_{B}b\longrightarrow ab_{0}\otimes b_{1}$. Therefore \cite{Masuoka03}
provides us a first example of a non-trivial Galois extension over a
coquasi-Hopf algebra.

\begin{example}
\rm%
(\cite{Albuquerque99a}) Let $G$ be any group and $\omega :G\times G\times
G\longrightarrow \Bbbk $ an invertible normalized cocycle. The the usual
group algebra $H=\Bbbk G$ becomes a coquasi-Hopf algebra by keeping the
ordinary operations, but with coassociator $\omega $ (linearly extended to $%
\Bbbk G^{\otimes 3}$) and linear maps $\alpha =\varepsilon $ and $\beta $
given by $\beta (g)=\omega ^{-1}(g,g^{-1},g)$, for any $g\in G$. As the
coalgebra structure is not modified, a $\Bbbk G$-coaction means precisely a $%
G$-graduation. Therefore, the notion of an $H$-comodule algebra becomes in
this case: a $G$-graded vector space $A=\oplus _{g\in G}A_{g}$, endowed with
a unit and a multiplication $"\cdot ":A\otimes A\rightarrow A$ such that $%
A_{g}A_{h}\subseteq A_{gh}$ for all $g,h\in G$, and associative in the sense
that%
\begin{equation*}
(a\cdot b)\cdot c=a\cdot (b\cdot c)\omega (\left\vert a\right\vert
,\left\vert b\right\vert ,\left\vert c\right\vert )
\end{equation*}%
for all homogeneous elements $a,b,c\in A$. The coinvariants $A^{coH}$ are
exactly $A_{e}$, where $e$ is the neutral element of $G$. We have then the
following:
\end{example}

\begin{proposition}
The extension $A_{e}\subseteq A$ is Galois (in the sense of Definition \ref%
{def Galois}) if and only if it is strongly graded.
\end{proposition}

\begin{proof}
Notice first that $A$ is strongly graded $\Longleftrightarrow $ $%
A_{g}A_{g^{-1}}=A_{e}$ for any $g\in G$. One inclusion is obvious, and for
the other we shall use the associativity rule of $A$: 
\begin{eqnarray*}
A_{gh} &\subseteq &A_{gh}A_{e}\subseteq
A_{gh}(A_{h^{-1}}A_{h})=(A_{gh}A_{h^{-1}})A_{h}\omega ^{-1}(gh,h^{-1},h) \\
&\subseteq &A_{g}A_{h}\omega ^{-1}(gh,h^{-1},h)\subseteq A_{g}A_{h}
\end{eqnarray*}

Now the proof follows as in the Hopf case.
\end{proof}

This result generalizes the Ulbrich's well-known example in the Hopf algebra
case (\cite{Ulbrich81}), and it is the first confirmation of the fact that
our definition of a Galois extension is the correct one.

\begin{example}
\rm%
Another example of Galois extension can be found in \cite{Balan07a}.
Starting from a coquasi-Hopf algebra $H$ and an associative algebra $R$
endowed with an $H$-weak action and a 2-cocycle $\sigma :H\otimes
H\longrightarrow R$, we can construct the crossed product $R\#_{\sigma }H$,
generalizing the Hopf case. This is a Galois extension of $R$ in the sense
of Definition \ref{def Galois}. Also, the Galois extension mentioned above
from \cite{Masuoka03} is precisely a particular case of our crossed product
construction.
\end{example}

Now, remember that a Hopf algebra $H$ can be seen as a right $H$-comodule
algebra via $\Delta $ and usual multiplication. The coinvariants are $\Bbbk
1_{H}\simeq \Bbbk $. Moreover, this extension is $H$-Galois (by \cite%
{Dascalescu00}) (actually, any bialgebra $H$ is a comodule algebra in this
way, and it is a Hopf algebra if and only if it is Galois). Now, working
with a coquasi-Hopf algebra $H$ still gives us a right comodule, but no
longer an algebra in the monoidal category $\mathcal{M}^{H}$ with the usual
multiplication. If we try to deform the multiplication on $H$ via a twist $%
\tau $ as in (\ref{multiplic Atau}), then $(H,\bullet _{\tau },\Delta )$ is
a right $H$-comodule algebra if and only if $\omega _{\tau }$ is trivial,
i.e. $H_{\tau }$ is a Hopf algebra. It is unclear to the author for the
moment which multiplication structure should be defined on $H$ such that we
get a right $H$-comodule algebra, which in the Hopf case should reduce to
ordinary multiplication. Moreover, this new multiplication should provide an
example of Galois extension $\Bbbk \subseteq H$.

\begin{remark}
\label{can op}%
\rm%
Let $H$ a coquasi-Hopf algebra with bijective antipode and $A$ a right $H$%
-comodule algebra. Notice that $\mathcal{M}^{H}\simeq {}^{H^{cop}}\mathcal{M}
$ as monoidal categories. Using also the isomorphism from (\ref{left-right}%
), it follows that $_{A}\mathcal{M}^{H}\simeq {}\mathcal{M}%
_{A^{op}}^{H^{op}} $ for $A$ a right $H$-comodule algebra. For completeness,
we remark that the corresponding Galois map for $H^{op}$-extension $%
B^{op}\subseteq A^{op}$ will be 
\begin{equation}
can^{\prime }(a\otimes _{B}b)=a_{0}b_{0}\otimes a_{4}\omega
(S^{-1}(a_{3})\beta S^{-1}(a_{2}),a_{1},b_{1})  \label{can prim}
\end{equation}%
Then we get:
\end{remark}

\begin{lemma}
\label{can bij implica canprim bij}The map $can^{\prime }$ is bijective if
and only if $can$ is bijective.
\end{lemma}

\begin{proof}
Consider the map $\Xi :A\otimes H\longrightarrow A\otimes H$, $\Xi (a\otimes
h)=a_{0}\otimes a_{3}S(h_{1})\omega ^{-1}(h_{3},a_{2}\beta (a_{3}),S(a_{4}))%
\mathbf{f}(h_{2},S^{-1}(a_{1}))$. Then $\Xi \circ can=can^{\prime }$. Also
one may check that $\Xi $ is bijective, with inverse given by 
\begin{equation*}
\Xi ^{-1}(a\otimes h)=a_{0}\otimes S^{-1}(h_{1})a_{5}\mathbf{f}%
^{-1}(S^{-1}(h_{2}),a_{4})\omega (a_{1}\beta (a_{2}),S(a_{3}),h_{3})
\end{equation*}%
Notice that in the case of a Hopf algebra the map $can^{\prime }$ reduces to
the usual formula $a\otimes _{B}b=a_{0}b\otimes _{B}a_{1}$.
\end{proof}

\begin{remark}
\rm%
\label{invarianta la twist}If $A$ is a right $H$-comodule algebra and $\tau $
a twist for $H$, we may consider the twisted comodule algebra $A_{\tau
^{-1}} $ as in Remark \ref{twist comodule algebra}. The comodule structure
being the same, $A$ and $A_{\tau ^{-1}}$ will have same coinvariants $B$
(but over different coquasi-Hopf algebras). Then we have the following:
\end{remark}

\begin{proposition}
The extension $B\subseteq A$ is $H$-Galois if and only if $B\subseteq
A_{\tau ^{-1}}$ is $H_{\tau }$-Galois.
\end{proposition}

\begin{proof}
The canonical Galois map for the extension $B\subseteq A_{\tau ^{-1}}$ is%
\begin{eqnarray}
can_{\tau }(a\otimes _{B}b) &=&a_{0}\cdot _{\tau }b_{0}\otimes \omega _{\tau
}^{-1}(a_{1},b_{1}\beta _{\tau }(b_{2}),S(b_{3}))b_{4}  \notag \\
(\ref{multiplic Atau}),(\ref{asociator Htau}),(\ref{beta Htau})
&=&a_{0}b_{0}\tau ^{-1}(a_{1},b_{1})\otimes \tau (a_{2},b_{2})\tau
(a_{3}b_{3},S(b_{13}))\omega ^{-1}(a_{4},b_{4},S(b_{12}))  \notag \\
&&\tau ^{-1}(a_{5},b_{5}S(b_{11}))\tau ^{-1}(b_{6},S(b_{10}))b_{14}\tau
(b_{7}\beta (b_{8}),S(b_{9}))  \notag \\
&=&a_{0}b_{0}\otimes \tau (a_{1}b_{1},S(b_{5}))\omega
^{-1}(a_{2},b_{2},\beta (b_{3})S(b_{4}))b_{6}  \label{can twist}
\end{eqnarray}%
for any $a,b\in A$. Consider now the linear map 
\begin{equation}
\vartheta :A\otimes H\longrightarrow A\otimes H,\vartheta (a\otimes
h)=a_{0}\otimes h_{2}\tau (a_{1},S(h_{1}))  \label{twist can}
\end{equation}%
It is easy to check that $\vartheta $ is bijective, with inverse $\vartheta
^{-1}(a\otimes h)=a_{0}\otimes h_{2}\tau ^{-1}(a_{1},S(h_{1}))$. Then the
following relation hold: $can_{\tau }=\vartheta \circ can$, which tells us
that both extensions will be simultaneously Galois.
\end{proof}

It follows from the previous Remark that if $H$ is a Hopf algebra and $%
B\subseteq A$ is a $H$-Galois extension in the classical sense, then for any
nontrivial twist $\tau \in (H\mathcal{\otimes }H)^{\ast }$, the extension $%
B\subseteq A_{\tau ^{-1}}$ will be $H_{\tau }$-Galois in the sense of our
Definition. Hence all known examples of Hopf-Galois extensions fit in our
picture.

\begin{example}
\rm%
Again, let $G$ a group and $\tau :G\times G\longrightarrow \Bbbk $ an
invertible normalized map. Then $H=\Bbbk G$ is a Hopf algebra and $A=\Bbbk G$
is an $H$-comodule algebra via comultiplication. Using the twist obtained by
extending $\tau $, it follows that $H_{\tau }$ is a coquasi-Hopf algebra and 
$A_{\tau ^{-1}}$ an $H_{\tau }$-comodule algebra. According to Remark \ref%
{invarianta la twist}, $A_{\tau ^{-1}}$ will be a Galois extension of $\Bbbk 
$. In particular, taking $G=(%
\mathbb{Z}
_{2})^{n}$, it follows that all Cayley algebras (as in \cite{Albuquerque99})
are Galois extensions over a coquasi-Hopf algebra.
\end{example}

\begin{remark}
\rm%
In the Hopf algebra case, the Galois map arises naturally as the evaluation
map $Hom_{A}^{H}(A,A\otimes H)\otimes _{B}A\longrightarrow A\otimes H$, from
the adjunction between the induced and the coinvariant functor, applied to
the relative Hopf module $A_{\bullet }\otimes H_{\bullet }^{\bullet }$. We
shall see that a similar result holds here, explaining thus the formula
chosen for $can$. We need first some work. For the beginning, it is not
obvious which $(H,A)$-Hopf module structure can be defined on $A\otimes H$
to generalize the one in the Hopf case. We shall assume that the antipode of 
$H$ is bijective, and obtaining the following:
\end{remark}

\begin{lemma}
\label{izom str urata de a-modul drept}The map 
\begin{eqnarray*}
\eta &:&H^{\bullet }\otimes A^{\bullet }\longrightarrow A\otimes H^{\bullet
}, \\
h\otimes a &\longrightarrow &a_{0}\otimes \omega (h_{1},a_{3},\alpha
S^{-1}(a_{2})S^{-1}(a_{1}))h_{2}a_{4}
\end{eqnarray*}%
is a right $H$-comodule isomorphism, where $H\otimes A$ is a comodule via
the codiagonal structure (i.e. $\rho _{H\otimes A}(h\otimes a)=h_{1}\otimes
a_{0}\otimes h_{2}a_{1}$) and $A\otimes H$ has the induced comodule
structure from the one of $H$.
\end{lemma}

\begin{proof}
The inverse for $\eta $ is given by $\eta ^{-1}(a\otimes
h)=h_{1}S^{-1}(a_{4})\otimes a_{0}\omega ^{-1}(h_{2},S^{-1}(a_{3})\beta
S^{-1}(a_{2}),a_{1})$.
\end{proof}

\begin{corollary}
\label{obs str urata de a-modul drept} Via the previous isomorphism, $%
A\otimes H$ becomes a right $A$-module in $\mathcal{M}^{H}$.
\end{corollary}

\begin{proof}
As $H$ is a right $H$-comodule via $\Delta $, $H^{\bullet }\otimes
A_{\bullet }^{\bullet }$ is naturally the right $A$-module induced in $%
\mathcal{M}^{H}$, with structures%
\begin{eqnarray*}
\rho _{H\otimes A}(h\otimes a) &=&h_{1}\otimes a_{0}\otimes h_{2}a_{1} \\
(h\otimes a)b &=&h_{1}\otimes a_{0}b_{0}\omega (h_{2},a_{1},b_{1})
\end{eqnarray*}%
for all $h\in H$, $a,b\in A$. Using $\eta $, the $A$-module structure can be
transferred on $A\otimes H$. Let's see how the multiplication formula with
elements of $A$ looks like:%
\begin{allowdisplaybreaks}%
\begin{eqnarray*}
(a\otimes h)\otimes b &\longrightarrow &[h_{1}S^{-1}(a_{4})\otimes
a_{0}\omega ^{-1}(h_{2},S^{-1}(a_{3})\beta S^{-1}(a_{2}),a_{1})]b \\
&=&h_{1}S^{-1}(a_{6})\otimes a_{0}b_{0}\omega
(h_{2}S^{-1}(a_{5}),a_{1},b_{1})\omega ^{-1}(h_{3},S^{-1}(a_{4})\beta
S^{-1}(a_{3}),a_{2}) \\
(\ref{cocycle omega}),(\ref{IdbetaS}) &=&h_{1}S^{-1}(a_{6})\otimes
a_{0}b_{0}\omega ^{-1}(h_{2},S^{-1}(a_{5}),a_{1}b_{1})\omega
(S^{-1}(a_{4}),a_{2},b_{2})\beta S^{-1}(a_{3}) \\
&\longrightarrow &a_{0}b_{0}\otimes \lbrack h_{2}S^{-1}(a_{10})](a_{4}b_{4})
\\
&&\omega (h_{1}S^{-1}(a_{11}),a_{3}b_{3},\alpha
S^{-1}(a_{2}b_{2})S^{-1}(a_{1}b_{1})) \\
&&\omega ^{-1}(h_{3},S^{-1}(a_{9}),a_{5}b_{5})\omega
(S^{-1}(a_{8}),a_{6},b_{6})\beta S^{-1}(a_{7}) \\
(\ref{asociat multipl}),(\ref{asociat multipl}),(\ref{IdbetaS})
&=&a_{0}b_{0}\otimes h_{3}b_{6}\omega (S^{-1}(a_{7}),a_{5},b_{5})\omega
^{-1}(h_{2},S^{-1}(a_{8}),a_{4}b_{4}) \\
&&\omega (h_{1}S^{-1}(a_{9}),a_{3}b_{3},S^{-1}(a_{1}b_{1}))\alpha
S^{-1}(a_{2}b_{2})\beta S^{-1}(a_{6}) \\
(\ref{asociat multipl}),(\ref{SalfaId}) &=&a_{0}b_{0}\otimes
h_{2}b_{7}\omega (S^{-1}(a_{10}),a_{4}b_{4},S^{-1}(a_{2}b_{2}))\omega
(h_{1},S^{-1}(a_{9})(a_{5}b_{5}), \\
&&S^{-1}(a_{1}b_{1}))\omega (S^{-1}(a_{8}),a_{6},b_{6})\alpha
S^{-1}(a_{3}b_{3})\beta S^{-1}(a_{7}) \\
(\ref{cocycle omega}),(\ref{IdbetaS}) &=&a_{0}b_{0}\otimes h_{2}b_{7}\omega
(S^{-1}(a_{8}),a_{4}b_{4},S^{-1}(a_{2}b_{2}))\omega
(S^{-1}(a_{7}),a_{5},b_{5}) \\
&&\omega (h_{1},b_{6},S^{-1}(a_{1}b_{1}))\alpha S^{-1}(a_{3}b_{3})\beta
S^{-1}(a_{6})
\end{eqnarray*}%
\end{allowdisplaybreaks}%

Therefore, we have 
\begin{eqnarray}
(a\otimes h)b &=&a_{0}b_{0}\otimes h_{2}b_{7}\omega
(S^{-1}(a_{8}),a_{4}b_{4},\alpha S^{-1}(a_{3}b_{3})S^{-1}(a_{2}b_{2}))\omega
(\beta S^{-1}(a_{6})S^{-1}(a_{7}),a_{5},b_{5})  \notag \\
&&\omega (h_{1},b_{6},S^{-1}(a_{1}b_{1}))  \label{str urata de a-modul drept}
\end{eqnarray}%
for any $h\in H$, $a,b\in A$. Seems to be complicated, but in the Hopf
algebra case it simply reduces to $(a\otimes h)b=ab_{0}\otimes hb_{1}$.
\end{proof}

We come back now to the counit of the adjunction, applied to the Hopf module 
$A\otimes H$. The coinvariants are $(A\otimes H)^{coH}=A\otimes \Bbbk
1_{H}\simeq A$, as the coaction takes place only on the second component and 
$\Bbbk $ is a commutative field. Hence 
\begin{allowdisplaybreaks}%
\begin{eqnarray}
\varepsilon _{A\otimes H} &:&A\otimes _{B}A\longrightarrow A\otimes H  \notag
\\
\varepsilon _{A\otimes H}(a\otimes _{B}b) &=&a_{0}b_{0}\otimes b_{5}\omega
(S^{-1}(a_{7}),a_{3}b_{3},\alpha S^{-1}(a_{2}b_{2})S^{-1}(a_{1}b_{1})) 
\notag \\
&&\omega (\beta S^{-1}(a_{5})S^{-1}(a_{6}),a_{4},b_{4})  \notag \\
(\ref{hdeltaS-1(a)h-1=S-1tensorS-1(delta cop(a))}) &=&a_{0}b_{0}\otimes
b_{7}\omega (S^{-1}(a_{9}),a_{5}b_{5},S^{-1}(b_{2})S^{-1}(a_{2}))\omega
(S^{-1}(a_{8}),a_{6},b_{6})  \notag \\
&&\alpha (S^{-1}(b_{3})S^{-1}(a_{3}))\beta S^{-1}(a_{7})\mathbf{f}%
^{(-1)}(S^{-1}(b_{1}),S^{-1}(a_{1}))  \notag \\
&&\mathbf{f}(S^{-1}(b_{4}),S^{-1}(a_{4}))  \notag \\
(\ref{f*alfa=gama, beta*f-1=delta}),(\ref{formula p}) &=&a_{0}b_{0}\otimes
b_{10}\mathbf{f}^{(-1)}(S^{-1}(b_{1}),S^{-1}(a_{1}))\beta
S^{-1}(a_{9})\alpha S^{-1}(b_{5})\alpha S^{-1}(a_{4})  \notag \\
&&\omega (a_{5},b_{6},S^{-1}(b_{4}))\omega
^{-1}(a_{6}b_{7},S^{-1}(b_{3}),S^{-1}(a_{3}))  \notag \\
&&\omega (S^{-1}(a_{11}),a_{7}b_{8},S^{-1}(b_{2})S^{-1}(a_{2}))\omega
(S^{-1}(a_{10}),a_{8},b_{9})  \notag \\
(\ref{cocycle omega}),(\ref{asociat multipl}),(\ref{SalfaId}),(\ref{IdbetaS}%
) &=&a_{0}b_{0}\otimes b_{10}\mathbf{f}^{(-1)}(S^{-1}(b_{1}),S^{-1}(a_{1}))%
\beta S^{-1}(a_{9})\alpha S^{-1}(b_{5})\alpha S^{-1}(a_{4})  \notag \\
&&\omega (S^{-1}(a_{12}),a_{5},S^{-1}(a_{3}))\omega
(a_{6},b_{6},S^{-1}(b_{4}))  \notag \\
&&\omega (S^{-1}(a_{11}),a_{7}b_{7},S^{-1}(b_{3}))\omega
(S^{-1}(a_{10}),a_{8},b_{8})  \notag \\
&&\omega ^{-1}(b_{9},S^{-1}(b_{2}),S^{-1}(a_{2}))  \notag \\
(\ref{cocycle omega}),(\ref{SalfaId}),(\ref{IdbetaS}),(\ref{omega anihileaza
S}) &=&a_{0}b_{0}\otimes b_{5}\mathbf{f}^{(-1)}(S^{-1}(b_{1}),S^{-1}(a_{1}))%
\alpha S^{-1}(b_{3})  \notag \\
&&\omega (b_{4},S^{-1}(b_{2}),S^{-1}(a_{2}))  \notag \\
(\ref{relatie UL pR h}) &=&a_{0}b_{0}\otimes b_{6}\mathbf{f}%
(b_{5},S^{-1}(a_{1}b_{1}))\omega ^{-1}(a_{2},b_{2}\beta (b_{3}),S(b_{4}))
\label{epsilon pt AtensorH}
\end{eqnarray}%
\end{allowdisplaybreaks}%

But according to (\ref{can twist}), this is precisely $can_{\widetilde{%
\mathbf{f}}}$, the Galois map twisted by $\widetilde{\mathbf{f}}$, where the
twist $\widetilde{\mathbf{f}}$ was introduced in relation (\ref{twist h}).
From the Remark \ref{invarianta la twist} it follows that:

\begin{corollary}
\label{counit bij implies Galois map bij}The bijectivity of $\varepsilon
_{A\otimes H}$ implies that $B\subseteq A$ is Galois.
\end{corollary}

\begin{remark}
\rm%
We could had used formula \ref{epsilon pt AtensorH} as an alternative Galois
map, but we preferred the formula from Definition \ref{def Galois} to avoid
the presence of the twist and simplify computations.
\end{remark}

We shall further need some properties of the Galois map, analogs to those in 
\cite{Schneider90}:

\begin{proposition}
The morphism $can$ satisfies the following:

(1) \label{prop can h colin}It is right $H$-colinear, where the right
comodule structure on both spaces is given from the second tensorand.

(2) \label{can trivial}For any $a\in A$, $can(1\otimes _{B}a)=a_{0}\otimes
\beta (a_{1})a_{2}$.

(3) \label{schneider}It is also right $H$-colinear, but with respect to the
following coactions: $\hat{\rho}(a\otimes _{B}b)=a_{0}\otimes _{B}b\otimes
a_{1}$, for $a\otimes _{B}b\in A\otimes _{B}A$, respectively $\check{\rho}%
(a\otimes h)=a_{0}\otimes h_{2}\otimes a_{1}S(h_{1})$, where $a\otimes h\in
A\otimes H$.

(4) \label{prop can a lin}If $can$ is bijective, then 
\begin{equation*}
(c\otimes 1_{H})can^{-1}(d\otimes h)=can^{-1}(c_{0}d_{0}\otimes h_{2})\omega
^{-1}(c_{1},d_{1},S(h_{1}))
\end{equation*}%
for any $c,d\in A$, $h\in H$.

(5) If the extension is Galois, denote $can^{-1}(1_{A}\otimes
h)=\sum\limits_{i}l_{i}(h)\otimes _{B}r_{i}(h)$. Then

(5.1) \label{prop can-1 h colin}$\sum\limits_{i}l_{i}(h_{1})\otimes
_{B}r_{i}(h_{1})\otimes h_{2}=\sum\limits_{i}l_{i}(h)\otimes
_{B}r_{i}(h)_{0}\otimes r_{i}(h)_{1}$.

(5.2) \label{l(h)r(h0=epsilon(h0}$\sum\limits_{i}l_{i}(h)r_{i}(h)=\alpha
(h)1_{A}$.

(5.3) \label{schneider can-1}$\sum\limits_{i}l_{i}(h)_{0}\otimes
_{B}r_{i}(h)\otimes l_{i}(h)_{1}=\sum\limits_{i}l_{i}(h_{2})\otimes
_{B}r_{i}(h_{2})\otimes S(h_{1})$.

(5.4) \label{prop can}$\sum\limits_{i}a_{0}\beta (a_{1})l_{i}(a_{2})\otimes
_{B}r_{i}(a_{2})=1_{A}\otimes _{B}a$.

(5.5) $\sum\limits_{i}l_{i}(hg)\otimes _{B}r_{i}(hg)=\sum\limits_{i,j}%
\mathbf{f}^{-1}(h_{1},g_{1})l_{i}(g_{2})l_{j}(h_{2})\otimes
_{B}r_{j}(h_{2})r_{i}(g_{2})$

for all $h,g\in H$, $a\in A$.

(6)The map $can^{\prime }$ from Remark \ref{can op} is also right $H$%
-colinear, where the right comodule structure on $A\otimes _{B}A$ is given
from the first tensorand, and $A\otimes H$ is a right $H$-comodule via $%
I_{A}\otimes \Delta $.
\end{proposition}

\begin{proof}
(1) We have 
\begin{eqnarray*}
(I\otimes \Delta )can(a\otimes _{B}b) &=&a_{0}b_{0}\otimes \omega
^{-1}(a_{1},b_{1}\beta (b_{2}),S(b_{3}))b_{4}\otimes b_{5} \\
&=&can(a\otimes _{B}b_{0})\otimes b_{1}
\end{eqnarray*}

(2) Obvious.

(3) Remark first that $\hat{\rho}$ and $\check{\rho}$ are indeed right $H$%
-comodule structures. Then, for any $a,b\in A$, we compute%
\begin{eqnarray*}
\check{\rho}can(a\otimes _{B}b) &=&\check{\rho}(a_{o}b_{0}\otimes \omega
^{-1}(a_{1},b_{1}\beta (b_{2}),S(b_{3})b_{4}) \\
&=&a_{0}b_{0}\otimes b_{6}\otimes (a_{1}b_{1})S(b_{5})\omega
^{-1}(a_{2},b_{2}\beta (b_{3}),S(b_{4})) \\
(\ref{asociat multipl}) &=&a_{0}b_{0}\otimes b_{6}\otimes
a_{2}(b_{2}S(b_{4}))\omega ^{-1}(a_{1},b_{1}\beta (b_{3}),S(b_{5})) \\
&=&a_{0}b_{0}\otimes b_{4}\otimes a_{2}\omega ^{-1}(a_{1},b_{1}\beta
(b_{2}),S(b_{3})) \\
&=&(can\otimes I_{H})(a_{0}\otimes _{B}b\otimes a_{1}) \\
&=&(can\otimes I_{H})\hat{\rho}(a\otimes _{B}b)
\end{eqnarray*}

(4) We get that%
\begin{eqnarray*}
can(ca\otimes _{B}b) &=&(c_{0}a_{0})b_{0}\otimes \omega
^{-1}(c_{1}a_{1},b_{1}\beta (b_{2}),S(b_{3}))b_{4} \\
&=&(c_{0}a_{0})b_{0}\otimes \omega ^{-1}(c_{1},a_{1},b_{1})\omega
^{-1}(c_{2},a_{2}b_{2},S(b_{6})_{1})\omega ^{-1}(a_{3},b_{3},S(b_{6})_{2}) \\
&&\omega (c_{3},a_{4},b_{4}S(b_{6})_{3})\beta (b_{5})b_{7} \\
&=&(c_{0}a_{0})b_{0}\otimes \omega ^{-1}(c_{1},a_{1},b_{1})\omega
^{-1}(c_{2},a_{2}b_{2},S(b_{8}))\omega ^{-1}(a_{3},b_{3},S(b_{7})) \\
&&\omega (c_{3},a_{4},b_{4}\beta (b_{5})S(b_{6}))b_{9} \\
&=&(c_{0}a_{0})b_{0}\otimes \omega ^{-1}(c_{1},a_{1},b_{1})\omega
^{-1}(c_{2},a_{2}b_{2},S(b_{6}))\omega ^{-1}(a_{3},b_{3},S(b_{5}))\beta
(b_{4})b_{7} \\
&=&c_{0}(a_{0}b_{0})\otimes \omega ^{-1}(c_{1},a_{1}b_{1},S(b_{5}))\omega
^{-1}(a_{2},b_{2},S(b_{4}))\beta (b_{3})b_{6}
\end{eqnarray*}%
for any $a,b,c\in A$. Now, if we denote $can^{-1}(d\otimes
h)=\sum\limits_{i}a_{i}\otimes _{B}b_{i}$, then%
\begin{equation*}
(\rho _{A}\otimes \Delta )(d\otimes h)=\sum_{i}a_{i0}b_{i0}\otimes
a_{i1}b_{i1}\otimes \omega ^{-1}(a_{i2},b_{i2}\beta
(b_{i3}),S(b_{i4}))b_{i5}\otimes b_{i6}
\end{equation*}%
Using this, we easily deduce that 
\begin{eqnarray*}
can(c\cdot can^{-1}(d\otimes h)) &=&can(\sum_{i}ca_{i}\otimes _{B}b_{i}) \\
&=&c_{0}(a_{i0}b_{i0})\otimes \omega
^{-1}(c_{1},a_{i1}b_{i1},S(b_{i5}))\omega
^{-1}(a_{i2},b_{i2},S(b_{i4}))\beta (b_{i3})b_{i6} \\
&=&c_{0}d_{0}\omega ^{-1}(c_{1},d_{1},S(h_{1}))\otimes h_{2}
\end{eqnarray*}

(5.1) It follows from (1).

(5.2) We shall check first the formula:%
\begin{eqnarray*}
\sum_{i}l_{i}(h)r_{i}(h) &=&(I_{A}\otimes \varepsilon
)(l_{i}(h)_{0}r_{i}(h)_{0}\otimes \omega
(l_{i}(h)_{1}r_{i}(h)_{1},S(r_{i}(h)_{6}),\alpha (r_{i}(h)_{7})r_{i}(h)_{8})
\\
&&\omega ^{-1}(l_{i}(h)_{2},r_{i}(h)_{2}\beta
(r_{i}(h)_{3}),S(r_{i}(h)_{4}))r_{i}(h)_{5})
\end{eqnarray*}%
The left hand side can be also written as:%
\begin{allowdisplaybreaks}%
\begin{align*}
& \sum_{i}(I_{A}\otimes \varepsilon )(l_{i}(h)_{0}r_{i}(h)_{0}\otimes \omega
(l_{i}(h)_{1}r_{i}(h)_{1},S(r_{i}(h)_{6}),\alpha (r_{i}(h)_{7})r_{i}(h)_{8})
\\
& \omega ^{-1}(l_{i}(h)_{2},r_{i}(h)_{2}\beta
(r_{i}(h)_{3}),S(r_{i}(h)_{4}))r_{i}(h)_{5} \\
& =\sum_{i}l_{i}(h)_{0}r_{i}(h)_{0}\omega
(l_{i}(h)_{1}r_{i}(h)_{1},S(r_{i}(h)_{5}),\alpha (r_{i}(h)_{6})r_{i}(h)_{7})
\\
& \omega ^{-1}(l_{i}(h)_{2},r_{i}(h)_{2}\beta (r_{i}(h)_{3}),S(r_{i}(h)_{4}))
\\
(\ref{cocycle omega})& =\sum_{i}l_{i}(h)_{0}r_{i}(h)_{0}\omega
^{-1}(l_{i}(h)_{1},r_{i}(h)_{1},S(r_{i}(h)_{5})_{1}r_{i}(h)_{7_{1}})\omega
(r_{i}(h)_{2},S(r_{i}(h)_{5})_{2},r_{i}(h)_{7_{2}}) \\
& \alpha (r_{i}(h)_{6})\beta (r_{i}(h)_{4})\omega
(l_{i}(h)_{2},r_{i}(h)_{3}S(r_{i}(h)_{5}),r_{i}(h)_{7_{3}}) \\
& =\sum_{i}l_{i}(h)r_{i}(h)_{0}\omega
(r_{i}(h)_{1},S(r_{i}(h)_{3}),r_{i}(h)_{5})\alpha (r_{i}(h)_{4})\beta
(r_{i}(h)_{2}) \\
& =\sum_{i}l_{i}(h)r_{i}(h)_{0}\varepsilon (r_{i}(h)_{1}) \\
& =\sum_{i}l_{i}(h)r_{i}(h)
\end{align*}%
\end{allowdisplaybreaks}%
But $\sum\limits_{i}l_{i}(h)_{0}r_{i}(h)_{0}\otimes \omega
^{-1}(l_{i}(h)_{1},r_{i}(h)_{1}\beta
(r_{i}(h)_{2}),S(r_{i}(h)_{3}))r_{i}(h)_{4}=1_{A}\otimes h$, therefore 
\begin{eqnarray*}
\sum_{i}l_{i}(h)r_{i}(h) &=&(I_{A}\otimes \varepsilon )(1_{A}\otimes
h_{1}\omega (1_{H},S(h_{2}),\alpha (h_{3})h_{4})) \\
&=&1_{A}\alpha (h)
\end{eqnarray*}

(5.3) It results from (3).

(5.4) We compute%
\begin{allowdisplaybreaks}
\begin{eqnarray*}
can(\sum_{i}a_{0}\beta (a_{1})l_{i}(a_{2})\otimes _{B}r_{i}(a_{2})) &=&\left[
a_{0}\beta (a_{2})l_{i}(a_{3})_{0}\right] r_{i}(a_{3})_{0}\otimes
r_{i}(a_{3})_{4} \\
&&\omega ^{-1}(a_{1}l_{i}(a_{3})_{1},r_{i}(a_{3})_{1}\beta
(r_{i}(a_{3})_{2}),S(r_{i}(a_{3})_{3}) \\
&=&a_{0}\left[ l_{i}(a_{4})_{0}r_{i}(a_{4})_{0}\right] \beta (a_{3})\omega
(a_{1},l_{i}(a_{4})_{1},r_{i}(a_{4})_{1})\otimes r_{i}(a_{4})_{5} \\
&&\omega ^{-1}(a_{2}l_{i}(a_{4})_{2},r_{i}(a_{4})_{2}\beta
(r_{i}(a_{4})_{3}),S(r_{i}(a_{4})_{4}) \\
&=&a_{0}\left[ l_{i}(a_{4})_{0}r_{i}(a_{4})_{0}\right] \beta (a_{3})\otimes
r_{i}(a_{4})_{8} \\
&&\omega
^{-1}(a_{1},l_{i}(a_{4})_{1}r_{i}(a_{4})_{1},S(r_{i}(a_{4})_{7}))\beta
(r_{i}(a_{4})_{4}) \\
&&\omega ^{-1}(l_{i}(a_{4})_{2},r_{i}(a_{4})_{2},S(r_{i}(a_{4})_{6})) \\
&&\omega (a_{2},l_{i}(a_{4})_{3},r_{i}(a_{4})_{3}S(r_{i}(a_{4})_{5}) \\
&=&a_{0}\left[ l_{i}(a_{3})_{0}r_{i}(a_{3})_{0}\right] \beta (a_{2})\otimes
r_{i}(a_{3})_{6} \\
&&\omega ^{-1}(a_{1},l_{i}(a_{3})_{1}r_{i}(a_{3})_{1},S(r_{i}(a_{3})_{5})) \\
&&\omega ^{-1}(l_{i}(a_{3})_{2},r_{i}(a_{3})_{2},S(r_{i}(a_{3})_{4}))\beta
(r_{i}(a_{3})_{3}) \\
&=&a_{0}\beta (a_{2})1_{A}\otimes a_{4}\omega ^{-1}(a_{1},1_{H},S(a_{3})) \\
&=&a_{0}\beta (a_{1})\otimes a_{2} \\
&=&can(1_{A}\otimes _{B}a)
\end{eqnarray*}%
\end{allowdisplaybreaks}%

(5.5) It is a consequence of the previous properties of $can$ and of the
properties (\ref{twist f})-(\ref{relatie UL pR h}) of the twist $\mathbf{f}$.

(6) Easy.
\end{proof}

\begin{remark}
\label{can morf de hopf mod}As a consequence of Proposition \ref{prop can a
lin}(3), and \ref{prop can a lin}(4), we obtain that $can$ is a morphism of
left Hopf modules, where $_{\bullet }A^{\bullet }\otimes _{B}A$ is an object
in $_{A}\mathcal{M}^{H}$ with structures given by the first tensorand, while 
$_{\bullet }A^{\bullet }\otimes H^{S,\bullet }$ is the induced module in $%
_{A}\mathcal{M}^{H}$. Here $H^{S,\bullet }$ is the right comodule structure
of $H$ deformed by the antipode $S$ (i.e. $_{\bullet }A^{\bullet }\otimes
H^{S,\bullet }$ is a left Hopf module with structure morphisms $a\otimes
h\longrightarrow a_{0}\otimes h_{2}\otimes a_{1}S(h_{1})$, $a(b\otimes
h)=a_{0}b_{0}\otimes h_{2}\omega ^{-1}(a_{1},b_{1},S(h_{1})$).
\end{remark}

\begin{definition}
Let $A$ a right $H$-comodule algebra and $\gamma :H\longrightarrow A$ a
colinear map. The extension $B\subseteq A$ is $(H,S)$-\textbf{cleft} with
respect to the cleaving map $\gamma $ if there is a linear map $\delta
_{\gamma ,S}:H\longrightarrow A$ such that 
\begin{eqnarray}
\rho (\delta _{\gamma ,S}(h)) &=&\delta _{\gamma ,S}(h_{2})\otimes S(h_{1})
\label{inversecleaving} \\
\delta _{\gamma ,S}(h_{1})\gamma (h_{2}) &=&\alpha (h)1_{A}
\label{convolutiedeltagama} \\
\gamma (h_{1})\beta (h_{2})\delta _{\gamma ,S}(h_{3}) &=&\varepsilon (h)1_{A}
\label{convolutiegamabetadelta}
\end{eqnarray}
\end{definition}

\begin{remark}
\rm%
\label{observatia cleft}(1) This definition of cleftness is slightly
different from the classical one. In the Hopf case, it is only required that 
$\gamma $ is convolution invertible (denote by $\delta $ the convolution
inverse of $\gamma $) and $H$-colinear. The property (\ref{inversecleaving})
appears naturally by passing from a bialgebra to a Hopf algebra.
Unfortunately, in our case the convolution product on $Hom(H,A)$ is no
longer associative, therefore a left inverse for $\gamma $ is not
necessarily a right inverse and the property (\ref{inversecleaving}) does
not seem to result from the other properties of $\gamma $. So we had to
state it separately.

(2) For a cleft comodule algebra $A$, the application $\delta _{\gamma ,S}$
depends on the antipode. If we change it to $(S^{\prime },\alpha ^{\prime
},\beta ^{\prime })$ as in (\ref{change antipode coquasi}) and define $%
\delta _{\gamma ,S^{\prime }}(h)=U(h_{1})\delta _{\gamma ,S}(h_{2})$, then
it follows immediately that $A$ is also $(H,S^{\prime })$-cleft. In the
sequel, we shall drop the subscripts for simplicity, considering the
antipode and the elements $\alpha $, $\beta $ fixed once for all.
\end{remark}

Recall that the "normal basis property" states that there is an isomorphism
of left $B$-modules, right $H$-comodules $A\simeq {}_{\bullet }B\otimes
H^{\bullet }$, where the dots are indicating the corresponding structures
for the tensor product. We shall keep the same definition for coquasi-Hopf
algebras, as nothing is changed.

\begin{theorem}
\label{doiTakBlattnerMontg}Let $H$ be a coquasi-Hopf algebra with bijective
antipode, $A$ a right $H$-comodule algebra with $B$ the subalgebra of
coinvariants. Then the following statements are equivalent:

(1) The extension $B\subseteq A$ is $H$-cleft;

(2) The Weak Structure Theorem holds and the extension has the normal basis
property;

(3) The extension $B\subseteq A$ is $H$-Galois and has the normal basis
property.

In this case, the categories $\mathcal{M}_{B}$ and $\mathcal{M}_{A}^{H}$ are
equivalent (the Strong Structure Theorem holds).
\end{theorem}

\begin{proof}
The proof of this theorem follows closely the original one for Hopf
algebras, due to Doi and Takeuchi (\cite{Doi86}), and Blattner and
Montgomery (\cite{Blattner89}), but we shall do the computations in detail,
because of their degree of difficulty.

(1) $\Longrightarrow $ (2) Define%
\begin{equation}
\nu :B\otimes H\longrightarrow A,\qquad \nu (b\otimes h)=b\gamma (h)
\label{normal basis izo}
\end{equation}%
It is obvious $B$-linear. As $\gamma $ is $H$-colinear, $\nu $ will also be.
We need an inverse for $\nu $. We take 
\begin{equation}
\nu ^{-1}(a)=a_{0}\delta (a_{1}\leftharpoonup \beta )\otimes a_{2}\text{, \ }%
a\in A  \label{normal basis inverse}
\end{equation}%
We need to show first that it is well-defined. For all $a\in A$, we have 
\begin{eqnarray}
\rho _{A}(a_{0}\delta (a_{1} &\leftharpoonup &\beta ))=a_{0_{0}}\delta
(a_{1}\leftharpoonup \beta )_{0}\otimes a_{0_{1}}\delta (a_{1}\leftharpoonup
\beta )_{1}  \notag \\
(\ref{inversecleaving}) &=&a_{0}\delta (a_{4})\otimes a_{1}\beta
(a_{2})S(a_{3})  \notag \\
&=&a_{0}\beta (a_{1})\delta (a_{2})\otimes 1_{H}  \notag \\
&=&a_{0}\delta (a_{1}\leftharpoonup \beta )\otimes 1_{H}
\label{proj to  coinvariants}
\end{eqnarray}%
meaning that $\func{Im}\nu ^{-1}\subseteq B\otimes H$. Let's check now that $%
\nu $ and $\nu ^{-1}$ are indeed inverses to each other: for all $a\in A$,
we compute%
\begin{eqnarray}
(\nu \circ \nu ^{-1})(a) &=&\nu (a_{0}\delta (a_{1}\leftharpoonup \beta
)\otimes a_{2})  \notag \\
&=&[a_{0}\delta (a_{1}\leftharpoonup \beta )]\gamma (a_{2})  \notag \\
&=&a_{0}[\delta (a_{2}\leftharpoonup \beta )_{0}\gamma (a_{3})_{0}]\omega
(a_{1},\delta (a_{2}\leftharpoonup \beta )_{1},\gamma (a_{3})_{1})  \notag \\
&=&a_{0}[\delta (a_{4})\gamma (a_{5})]\omega (a_{1},S(a_{3}),a_{6})\beta
(a_{2})  \notag \\
(\ref{convolutiedeltagama}) &=&a_{0}\omega (a_{1},S(a_{3}),a_{5})\beta
(a_{2})\alpha (a_{4})  \notag \\
&=&a  \label{normal basis prop is  bijective}
\end{eqnarray}%
Conversely, for $b\in B$ and $h\in H$ we get 
\begin{eqnarray*}
\nu ^{-1}\circ \nu (b\otimes h) &=&\nu ^{-1}(b\gamma (h)) \\
&=&b\gamma (h)_{0}\delta (\left( \gamma (h)_{1}\right) \leftharpoonup \beta
)\otimes \gamma (h)_{2} \\
&=&b\gamma (h_{1})\beta (h_{2})\delta (h_{3})\otimes h_{4} \\
(\ref{convolutiegamabetadelta}) &=&b\otimes h
\end{eqnarray*}%
Hence the extension $B\subseteq A$ has the normal basis property. It remains
only to show the bijectivity of the adjunction counit from Proposition \ref%
{functorul indus}. For a Hopf module $M\in \mathcal{M}_{A}^{H}$, define the
map $t_{M}:M\longrightarrow M$ by $t_{M}(m)=m_{0}\delta (m_{1}\leftharpoonup
\beta )$. Then we can see as in (\ref{proj to coinvariants}) that the image
of $t_{M}$ is in $M^{coH}$. Define now $\chi :M\longrightarrow
M^{coH}\otimes _{B}A$, $\chi (m)=t_{M}(m_{0})\otimes _{B}\gamma (m_{1})$.
Computing as in (\ref{normal basis prop is bijective}), we get that $\chi $
is an inverse for $\varepsilon _{M}$.

(2) $\Longrightarrow $ (3) It follows from Corollary (\ref{counit bij
implies Galois map bij}).

(3) $\Longrightarrow $ (1) Let $\nu :B\otimes H\longrightarrow A$ be the
isomorphism given by the normal basis property. Define $\gamma (h)=\nu
(1_{A}\otimes h)$. As $\nu $ is $H$-colinear, $\gamma $ will also be.

In order to get the second map $\delta $, we need some work first. Consider
the map $\Gamma =(I_{A}\otimes \varepsilon )\nu ^{-1}:A\longrightarrow B$.
Then $\Gamma $ is left $B$-linear, as $\nu ^{-1}$ is $B$-linear, and%
\begin{eqnarray}
\Gamma \gamma (h) &=&(I_{A}\otimes \varepsilon )\nu ^{-1}\nu (1_{A}\otimes h)
\notag \\
&=&\varepsilon (h)1_{A}  \label{g(gama)}
\end{eqnarray}

Now we may take $\delta (h)=m_{A}(I_{A}\otimes _{B}\Gamma )can^{-1}(1\otimes
h)=\sum\limits_{i}l_{i}(h)\Gamma (r_{i}(h))$, where $m_{A}$ is the
multiplication on $A$. We may then compute%
\begin{eqnarray*}
\gamma (h_{1})\beta (h_{2})\delta (h_{3}) &=&\underset{\in A}{\underbrace{%
\gamma (h_{1})}}m_{A}(I_{A}\otimes _{B}\Gamma )can^{-1}(1_{A}\otimes \beta
(h_{2})h_{3}) \\
&=&m_{A}(I_{A}\otimes _{B}\Gamma )[\gamma (h_{1})can^{-1}(1_{A}\otimes \beta
(h_{2})h_{3})] \\
(\text{Proposition }\ref{prop can a lin}) &=&m_{A}(I_{A}\otimes _{B}\Gamma
)can^{-1}(\gamma (h_{1})_{0}\otimes \beta (h_{2})h_{4})\omega ^{-1}(\gamma
(h_{1})_{1},1,S(h_{3})) \\
&=&m_{A}(I_{A}\otimes _{B}\Gamma )can^{-1}(\gamma (h_{1})\otimes \beta
(h_{2})h_{3}) \\
(\gamma \text{ is colinear}) &=&m_{A}(I_{A}\otimes _{B}\Gamma
)can^{-1}(\gamma (h)_{0}\otimes \beta (\gamma (h)_{1})\gamma (h)_{2}) \\
(\text{Proposition }\ref{can trivial}) &=&m_{A}(I_{A}\otimes _{B}\Gamma
)(1_{A}\otimes _{B}\gamma (h)) \\
&=&\Gamma \gamma (h) \\
&=&\varepsilon (h)1_{A}
\end{eqnarray*}

For the last formula, notice first that $H$-colinearity of $\nu $ implies%
\begin{eqnarray}
\nu ^{-1} &=&(I_{B}\otimes \varepsilon \otimes I_{H})(I_{B}\otimes \Delta
)\nu ^{-1}  \notag \\
&=&(I_{B}\otimes \varepsilon \otimes I_{H})(\nu ^{-1}\otimes I_{H})\rho _{A}
\notag \\
&=&(\Gamma \otimes I_{H})\rho _{A}  \label{fi-1=g  tensor I ro}
\end{eqnarray}

Now we may compute%
\begin{allowdisplaybreaks}%
\begin{eqnarray*}
\delta (h_{1})\gamma (h_{2}) &=&[m_{A}(I_{A}\otimes _{B}\Gamma
)can^{-1}(1_{A}\otimes h_{1})]\nu (1_{A}\otimes h_{2}) \\
&=&\sum_{i}[l_{i}(h_{1})\underset{\in B}{\underbrace{\Gamma (r_{i}(h_{1}))}}%
]\nu (1_{A}\otimes h_{2}) \\
(\text{Proposition }\ref{prop can-1 h colin}) &=&\sum_{i}[l_{i}(h)\underset{%
\in B}{\underbrace{\Gamma (r_{i}(h)_{0})}}]\nu (1_{A}\otimes r_{i}(h)_{1}) \\
&=&\sum_{i}l_{i}(h)[\Gamma (r_{i}(h)_{0})\nu (1_{A}\otimes r_{i}(h)_{1})] \\
&=&\sum_{i}l_{i}(h)\nu (\Gamma (r_{i}(h)_{0}\otimes r_{i}(h)_{1}) \\
(\ref{fi-1=g tensor I ro}) &=&\sum_{i}l_{i}(h)r_{i}(h) \\
(\text{Proposition }\ref{l(h)r(h0=epsilon(h0}) &=&\alpha (h)1_{A}
\end{eqnarray*}%
\end{allowdisplaybreaks}%
for all $h\in H$.

Finally, for any $h\in H$ we have 
\begin{eqnarray*}
\rho _{A}\delta (h) &=&\rho _{A}m_{A}(I_{A}\otimes _{B}\Gamma
)can^{-1}(1_{A}\otimes h) \\
&=&\rho _{A}(\sum_{i}l_{i}(h)\Gamma (r_{i}(h))) \\
&=&\sum_{i}l_{i}(h)_{0}\Gamma (r_{i}(h))\otimes l_{i}(h)_{1} \\
(\text{Proposition }\ref{schneider can-1}) &=&\sum_{i}l_{i}(h_{2})\Gamma
(r_{i}(h_{2}))\otimes S(h_{1}) \\
&=&\delta (h_{2})\otimes S(h_{1})
\end{eqnarray*}

For the remaining of the theorem, the proof is the same as in \cite{Doi86},
so we omit it.
\end{proof}

We shall prove now an imprimitivity statement which is the analogue of Doi's
and Takeuchi's theorem (\cite{Doi89}) and Schneider's theorem (\cite%
{Schneider90a}) for coquasi-Hopf algebras.

\begin{theorem}
\label{a galois + fidel plat echiv strong str th}Let $H$ be a coquasi-Hopf
algebra with bijective antipode, $A$ a right $H$-comodule algebra with $B$
the algebra of coinvariants. Then the following are equivalent:

(1) $A$ is faithfully flat as a left $B$-module, and $A$ is a Galois
extension of $B$.

(2) The functor of coinvariants and the induction functor are a pair of
inverse equivalences between $\mathcal{M}_{A}^{H}$ and $\mathcal{M}_{B}$.
\end{theorem}

\begin{proof}
(1) $\Longrightarrow $ (2) We need first a Lemma:

\begin{lemma}
\label{canM bijectiv}Let $H$ be a coquasi-Hopf algebra and $A$ a right
comodule algebra. Then the Galois map $can$ induces a natural right colinear
morphism $can_{M}:M\otimes _{B}A^{\bullet }\longrightarrow M\otimes
H^{\bullet }$, $can_{M}(m\otimes _{B}a)=m_{0}a_{0}\otimes \omega
^{-1}(m_{1},a_{1}\beta (a_{2}),S(a_{3}))a_{4}$. If $can$ is bijective, then $%
can_{M}$ is also bijective.
\end{lemma}

\begin{proof}
(of the \textbf{Lemma}) The definition of $can_{M}$ allows us to easily
check its colinearity. The hard part is the proof of the naturality and of
the bijectivity of $can_{M}$. This is not obvious at all, because $A$ is no
longer an associative algebra and the classical argument (i.e. tensoring
over $A$) is not working. In this case it is more appropriate to work in the
monoidal category of comodules. We refer to \cite{Ardizzoni07}\ for details
about tensor product over an algebra in a monoidal category. As $A$ is an
algebra in the monoidal abelian category $\mathcal{M}^{H}$, we may form the
tensor product $M\bigcirc _{A}M^{\prime }$ for any $M\in \mathcal{M}_{A}^{H}$%
, $M^{\prime }\in {}_{A}\mathcal{M}^{H}$ as the following equalizer 
\begin{equation*}
(M\otimes A)\otimes M^{\prime }\overset{\mu _{M}\otimes I_{M^{\prime }}}{%
\underset{(I_{M}\otimes \mu _{M^{\prime }})\phi _{M,A,M^{\prime }}}{%
\rightrightarrows }}M\otimes M^{\prime }\longrightarrow M\bigcirc
_{A}M^{\prime }\longrightarrow 0
\end{equation*}%
where $\mu _{M}$ and $\mu _{M^{\prime }}$ are the $A$-module structure maps
and $\phi $ is the coassociator.

We need now two particular left Hopf modules. One of them is $_{\bullet
}A^{\bullet }\otimes _{B}A$, with right $H$-coaction and left $A$-action on
the first component. For the other one, notice first that $S$ is a coalgebra
map. Therefore we may take $H$ as an object in $\mathcal{M}^{H}$ with $%
h\longrightarrow h_{2}\otimes S(h_{1})$, denoted $H^{S}$. Then we get an
induced left Hopf module $_{\bullet }A^{\bullet }\otimes H^{S,\bullet }\in $ 
$_{A}\mathcal{M}^{H}$, with structure maps 
\begin{eqnarray}
a\otimes h &\longrightarrow &a_{0}\otimes h_{2}\otimes a_{1}S(h_{1})
\label{comodul ob indus} \\
a(b\otimes h) &=&a_{0}b_{0}\otimes h_{2}\omega ^{-1}(a_{1},b_{1},S(h_{1}))
\label{A modul ob indus}
\end{eqnarray}

We can construct now the following diagram for any $M\in \mathcal{M}_{A}^{H}$%
:%
\begin{equation*}
\begin{array}{ccccc}
&  &  &  & M\otimes _{B}A \\ 
&  &  & F_{1}\nearrow \swarrow G_{1} & \widetilde{F}_{1}\uparrow \downarrow 
\widetilde{F}_{1}^{-1} \\ 
(M\otimes A)\otimes (A\otimes _{B}A) & \overset{\mu _{M}\otimes
(I_{A}\otimes _{B}I_{A})}{\underset{(I_{M}\otimes (m_{A}\otimes
_{B}I_{A})\phi _{M,A,A\otimes _{B}A}}{\rightrightarrows }} & M\otimes
(A\otimes _{B}A) & \overset{\pi _{1}}{\longrightarrow } & M\bigcirc
_{A}(A\otimes _{B}A)\longrightarrow 0 \\ 
&  & \downarrow _{I_{M}\otimes can} &  & \downarrow _{I_{M}\bigcirc _{A}can}
\\ 
(M\otimes A)\otimes (A\otimes H) & \overset{\mu _{M}\otimes (I_{A}\otimes
_{B}I_{A})}{\underset{(I_{M}\otimes \mu _{A\otimes H})\phi _{M,A,A\otimes H}}%
{\rightrightarrows }} & M\otimes (A\otimes H) & \overset{\pi _{2}}{%
\longrightarrow } & M\bigcirc _{A}(A\otimes H)\longrightarrow 0 \\ 
&  &  & F_{2}\searrow \nwarrow G_{2} & \widetilde{F}_{2}\downarrow \uparrow 
\widetilde{G}_{2} \\ 
&  &  &  & M\otimes H%
\end{array}%
\end{equation*}

The two rows are exact by definition of $\bigcirc _{A}$. The application $%
F_{1}:M_{\bullet }^{\bullet }\otimes (_{\bullet }A^{\bullet }\otimes
_{B}A)\longrightarrow M^{\bullet }\otimes _{B}A$, $F_{1}(m\otimes (a\otimes
_{B}b))=ma\otimes _{B}b$ is well-defined, right $H$-colinear and 
\begin{equation*}
F_{1}(\mu _{M}\otimes (I_{A}\otimes _{B}I_{A}))=F_{1}(I_{M}\otimes
(m_{A}\otimes _{B}I_{A}))\phi _{M,A,A\otimes _{B}A}
\end{equation*}%
(here $M^{\bullet }\otimes _{B}A$ is a right comodule via $\rho _{M}\otimes
_{B}I_{A}$, while $M_{\bullet }^{\bullet }\otimes (_{\bullet }A^{\bullet
}\otimes _{B}A)$ has the codiagonal comodule structure). Hence there is a
right $H$-comodule morphism $\widetilde{F}_{1}:M\bigcirc _{A}(A\otimes
_{B}A)\longrightarrow M\otimes _{B}A$ such that $\widetilde{F}_{1}\pi
_{1}=F_{1}$. Moreover, $F_{1}$ is an isomorphism with inverse $\widetilde{F}%
_{1}^{-1}=\pi _{1}G_{1}$, where $G_{1}(m\otimes _{B}a)=m\otimes
(1_{A}\otimes _{B}a)$.

As $_{\bullet }A^{\bullet }\otimes H^{S,\bullet }$ is an $A$-module induced
in $\mathcal{M}^{H}$, the colinear map $F_{2}:M_{\bullet }^{\bullet }\otimes
(_{\bullet }A^{\bullet }\otimes H^{S,\bullet })\longrightarrow M^{\bullet
}\otimes H^{S,\bullet }$, $F_{2}(m\otimes (a\otimes h))=m_{0}a_{0}\otimes
h_{2}\omega ^{-1}(m_{1},a_{1},S(h_{1}))$ factors through an isomorphism of
right $H$-comodules $\widetilde{F}_{2}:M\bigcirc _{A}(A\otimes H)\simeq
M\otimes H$ , with inverse $\widetilde{F}_{2}^{-1}=\pi _{2}G_{2}$, where $%
G_{2}(m\otimes h)=m\otimes (1_{A}\otimes h)$ (on $M^{\bullet }\otimes
H^{S,\bullet }$ we have again the codiagonal comodule structure).

According to Remark \ref{can morf de hopf mod}, the map $can:{}_{\bullet
}A^{\bullet }\otimes _{B}A\longrightarrow {}_{\bullet }A^{\bullet }\otimes
H^{S,\bullet }$ is a morphism in $_{A}\mathcal{M}^{H}$. Then $I_{M}\otimes
can$ induces a colinear map $I_{M}\bigcirc _{A}can:M\bigcirc _{A}(A\otimes
_{B}A)\longrightarrow M\bigcirc _{A}(A\otimes _{B}A)$.

Composing now the morphisms from the last column in the above diagram, we
obtain a natural map $can_{M}:M\otimes _{B}A\longrightarrow M\otimes H$,
which can be written as $can_{M}=\widetilde{F}_{2}(I_{M}\bigcirc _{A}can)%
\widetilde{F}_{1}^{-1}=F_{2}(I_{M}\otimes can)G_{1}$. This implies 
\begin{eqnarray*}
can_{M}(m\otimes _{B}a) &=&F_{2}(I_{M}\otimes can)G_{1}(m\otimes _{B}a) \\
&=&F_{2}(I_{M}\otimes can)(m\otimes (1_{A}\otimes _{B}a)) \\
&=&F_{2}(m\otimes (a_{0}\otimes \beta (a_{1})a_{2})) \\
&=&m_{0}a_{0}\otimes \omega ^{-1}(m_{1},a_{1}\beta (a_{2}),S(a_{3}))a_{4}
\end{eqnarray*}%
The last part of the Lemma is now obvious.
\end{proof}

(Proof of the\textbf{\ Theorem}) As in Lemma \ref{izom str urata de a-modul
drept}, we are able to show that for any $M\in \mathcal{M}_{A}^{H}$, $%
M\otimes H$ becomes an object in $\mathcal{M}_{A}^{H}$, with structure
morphisms as in Remark \ref{obs str urata de a-modul drept} (replacing the
elements of $A$ with elements of $M$). For later use, we write down
explicitly the used isomorphisms:%
\begin{eqnarray*}
&&H^{\bullet }\otimes M^{\bullet }\overset{\eta _{M}}{\underset{\eta
_{M}^{-1}}{\rightleftarrows }}M\otimes H^{\bullet }, \\
\eta _{M}(h\otimes m) &=&m_{0}\otimes \omega (h_{1},m_{3},\alpha
S^{-1}(m_{2})S^{-1}(m_{1}))h_{2}m_{4} \\
\eta _{M}^{-1}(m\otimes h) &=&h_{1}S^{-1}(m_{4})\otimes m_{0}\omega
^{-1}(h_{2},S^{-1}(m_{3})\beta S^{-1}(m_{2}),m_{1})
\end{eqnarray*}%
It follows that $\varepsilon _{M\otimes H}=(can_{M})_{\widetilde{\mathbf{f}}%
} $ is bijective. Repeating the argument with $M\otimes H$ instead of $M$,
we obtain the bijectivity of $\varepsilon _{(M\otimes H)\otimes
H}=(can_{M\otimes H})_{\widetilde{\mathbf{f}}}$. Now the trick is to put
this three maps together in a commutative diagram using their naturality: 
\begin{equation}
\begin{array}{cccccc}
0\longrightarrow & M^{coH}\otimes _{B}A & \longrightarrow & M\otimes _{B}A & 
\underset{}{\overset{}{\rightrightarrows }} & ((M\otimes H)\otimes _{B}A \\ 
& \varepsilon _{M}\downarrow &  & \varepsilon _{M\otimes H}\downarrow &  & 
\varepsilon _{(M\otimes H)\otimes H}\downarrow \\ 
0\longrightarrow & M & \overset{\widetilde{\rho }_{M}}{\longrightarrow } & 
M\otimes H & \underset{\widetilde{\rho }_{M}\otimes I_{H}}{\overset{%
I_{M}\otimes \widetilde{\Delta }}{\rightrightarrows }} & (M\otimes H)\otimes
H%
\end{array}
\label{diagrama counit}
\end{equation}%
On the bottom row, $\widetilde{\rho }_{M}(m)=m_{0}\otimes \alpha
S^{-1}(m_{1})m_{2}$ and $\widetilde{\Delta }(h)=h_{1}\otimes \alpha
S^{-1}(h_{2})h_{3}$ for $m\in M$, $h\in H$. The upper row contains their
images via the functor $(-)^{coH}\otimes _{B}A$. We need to see that $%
\widetilde{\rho }_{M}$, $\widetilde{\rho }_{M}\otimes I_{H}$ and $%
I_{M}\otimes \widetilde{\Delta }$ are morphisms of Hopf modules and the
bottom row is exact. It is easy to check the colinearity of these maps, as
the comodule structure on the tensor products $M\otimes H$ and $M\otimes
H\otimes H$ comes from the last tensorand. The difficult part is the right $%
A $-linearity, because of the unpleasant formula (\ref{str urata de a-modul
drept}) (applied to $M$, respectively $M\otimes H$). Instead of checking it
directly, we shall use the isomorphism which leads to the mentioned formula.
For $\widetilde{\rho }_{M}$, we have:%
\begin{equation*}
\begin{array}{ccccc}
M & \overset{\widetilde{\rho }_{M}}{\longrightarrow } & M\otimes H & \overset%
{\eta _{M}^{-1}}{\simeq } & H\otimes M%
\end{array}%
\end{equation*}%
and%
\begin{eqnarray}
m &\longrightarrow &m_{0}\otimes \alpha S^{-1}(m_{1})m_{2}\overset{\eta
_{M}^{-1}}{\longrightarrow }m_{6}\alpha S^{-1}(m_{5})S^{-1}(m_{4})\otimes
m_{0}\omega ^{-1}(m_{7},S^{-1}(m_{3})\beta S^{-1}(m_{2}),m_{1})  \notag \\
(\text{\ref{SalfaId}}) &=&\alpha S^{-1}(m_{4})1_{H}\otimes m_{0}\omega
^{-1}(m_{5},S^{-1}(m_{3})\beta S^{-1}(m_{2}),m_{1})  \notag \\
(\text{\ref{omega anihileaza S}}) &=&1_{H}\otimes m  \label{miu compus cu ro}
\end{eqnarray}%
Remembering that the $A$-module structure of $H\otimes M$ is the one induced
by $M$ (i.e. $(h\otimes m)a=h_{1}\otimes m_{0}a_{0}\omega
(h_{2},m_{1},a_{1}) $), it is now easy to verify the $A$-linearity for the
composed map $\eta _{M}^{-1}\widetilde{\rho }_{M}$.

For $\widetilde{\rho }_{M}\otimes I_{H}$, we need to compose three times
with the following isomorphisms%
\begin{equation*}
H^{\bullet }\otimes M_{\bullet }^{\bullet }\overset{\eta _{M}}{\simeq }%
M_{\bullet }\otimes H_{\bullet }^{\bullet }\overset{\widetilde{\rho }%
_{M}\otimes I_{H}}{\longrightarrow }(M_{\bullet }\otimes H_{\bullet
}^{\bullet })_{\bullet }\otimes H_{\bullet }^{\bullet }\overset{\eta
_{M\otimes H}}{\simeq }H^{\bullet }\otimes (M_{\bullet }\otimes H_{\bullet
}^{\bullet })_{\bullet }^{\bullet }\overset{I_{H}\otimes \eta _{M}^{-1}}{%
\simeq }H^{\bullet }\otimes (H^{\bullet }\otimes M_{\bullet }^{\bullet
})_{\bullet }^{\bullet }
\end{equation*}%
We use dots to indicate the structures (although the right $A$-structures
are not the classical ones, we find the notation more suggestive). We obtain 
\begin{eqnarray*}
h\otimes m &\longrightarrow &m_{0}\otimes \omega (h_{1},m_{3},\alpha
S^{-1}(m_{2})S^{-1}(m_{1}))h_{2}m_{4} \\
&\longrightarrow &m_{0}\otimes \alpha S^{-1}(m_{1})m_{2}\otimes \omega
(h_{1},m_{5},\alpha S^{-1}(m_{4})S^{-1}(m_{3}))h_{2}m_{6} \\
&\longrightarrow &h_{1}\otimes m_{0}\otimes \alpha S^{-1}(m_{1})m_{2}\omega
^{-1}(h_{3}m_{9},S^{-1}(m_{5})\beta S^{-1}(m_{4}),m_{3})\omega
(h_{2},m_{8}\alpha S^{-1}(m_{7}),S^{-1}(m_{6})) \\
(\text{\ref{cocycle omega}}) &=&h_{1}\otimes m_{0}\otimes \alpha
S^{-1}(m_{1})m_{2}\omega ^{-1}(h_{2},m_{11}S^{-1}(m_{9}),m_{3})\omega
^{-1}(m_{12},S^{-1}(m_{8}),m_{4}) \\
&&\omega (h_{3},m_{10},S^{-1}(m_{7})m_{5})\beta S^{-1}(m_{6})\alpha
S^{-1}(m_{10}) \\
(\text{\ref{SalfaId}, \ref{IdbetaS}}) &=&h\otimes m_{0}\otimes \alpha
S^{-1}(m_{1})m_{2}\omega ^{-1}(m_{7},S^{-1}(m_{5}),m_{3})\beta
S^{-1}(m_{4})\alpha S^{-1}(m_{6}) \\
(\text{\ref{omega anihileaza S}}) &=&h\otimes m_{0}\otimes \alpha
S^{-1}(m_{1})m_{2} \\
(\text{\ref{miu compus cu ro}}) &\longrightarrow &h\otimes 1_{H}\otimes m
\end{eqnarray*}%
and this is again right $A$-linear.

Finally, we repeat the above composition with $I_{M}\otimes \widetilde{%
\Delta }$ instead of $\widetilde{\rho }_{M}\otimes I_{H}$. We obtain 
\begin{eqnarray*}
h\otimes m &\longrightarrow &m_{0}\otimes \omega (h_{1},m_{3},\alpha
S^{-1}(m_{2})S^{-1}(m_{1}))h_{2}m_{4} \\
&\longrightarrow &m_{0}\otimes \omega (h_{1},m_{3},\alpha
S^{-1}(m_{2})S^{-1}(m_{1}))h_{2}m_{4}\otimes \alpha
S^{-1}(h_{3}m_{5})h_{4}m_{6} \\
&\longrightarrow &h_{8}m_{10}\alpha
S^{-1}(h_{7}m_{9})S^{-1}(h_{6}m_{8})\otimes m_{0}\otimes \omega
(h_{1},m_{3},\alpha S^{-1}(m_{2})S^{-1}(m_{1}))h_{2}m_{4} \\
&&\omega ^{-1}(h_{9}m_{11},S^{-1}(h_{5}m_{7})\beta
S^{-1}(h_{4}m_{6}),h_{3}m_{5}) \\
(\text{\ref{SalfaId}}) &=&1_{H}\otimes m_{0}\otimes \omega
(h_{1},m_{3},\alpha S^{-1}(m_{2})S^{-1}(m_{1}))h_{2}m_{4}\omega
^{-1}(h_{7}m_{9},S^{-1}(h_{5}m_{7}) \\
&&\beta S^{-1}(h_{4}m_{6}),h_{3}m_{5})\alpha S^{-1}(h_{6}m_{8}) \\
(\text{\ref{omega anihileaza S}}) &=&1_{H}\otimes m_{0}\otimes \omega
(h_{1},m_{3},\alpha S^{-1}(m_{2})S^{-1}(m_{1}))h_{2}m_{4} \\
&\longrightarrow &1_{H}\otimes h\otimes m
\end{eqnarray*}%
which respects the multiplication with elements of $A$.

We have to show now the exactness of the sequence. For the injectivity of $%
\widetilde{\rho }_{M}$: take $m\in M$ such that $m_{0}\otimes \alpha
S^{-1}(m_{1})m_{2}=0$. Now apply $\rho _{M}\otimes \lbrack (I_{H}\otimes
\Delta )\Delta ]$ to get 
\begin{equation*}
0=m_{0}\otimes m_{1}\otimes \alpha S^{-1}(m_{2})m_{3}\otimes m_{4}\otimes
m_{5}
\end{equation*}%
Finally, act on this by $I_{M}\otimes \omega ^{321}(S^{-1}\otimes
I_{H}\otimes \beta S^{-1}\otimes S^{-1})$. By (\ref{omega anihileaza S}) it
follows $0=m_{0}\varepsilon (m_{1})=m$.

Let's check now the exactness in $M\otimes H$. It is straightforward to see
that $(I_{M}\otimes \widetilde{\Delta })\widetilde{\rho }_{M}=(\widetilde{%
\rho }_{M}\otimes I_{H})\widetilde{\rho }_{M}$. Conversely, let $%
\sum_{i}m_{i}$ $\otimes h_{i}\in M\otimes H$ such that $\sum_{i}m_{i}\otimes
h_{i1}\otimes \alpha S^{-1}(h_{i2})h_{i3}=\sum_{i}m_{i0}\otimes \alpha
S^{-1}(m_{i1})m_{i2}\otimes h_{i}$. Apply $\Delta $ and $(I_{H}\otimes
\Delta )\Delta $ on the second, respectively last component, then act by $%
I_{M}\otimes I_{H}\otimes \omega ^{321}(S^{-1}\otimes I_{H}\otimes \beta
S^{-1}\otimes S^{-1})$. Again using (\ref{omega anihileaza S}) we obtain 
\begin{eqnarray*}
\sum_{i}m_{i}\otimes h_{i} &=&\sum_{i}m_{i0}\otimes \alpha
S^{-1}(m_{i1})m_{i2}\omega (S^{-1}(h_{i3})\beta
S^{-1}(h_{i2}),h_{i1},S^{-1}(m_{i3})) \\
&=&\widetilde{\rho }_{M}(\sum_{i}m_{i0}\omega (S^{-1}(h_{i3})\beta
S^{-1}(h_{i2}),h_{i1},S^{-1}(m_{i1}))
\end{eqnarray*}%
Therefore the diagram (\ref{diagrama counit})\ is commutative by the
naturality of $\varepsilon _{(-)}$. The upper row is exact because $%
(-)^{coH} $ is exact and $A$ is a faithfully flat $B$-module. As explained
above, $\varepsilon _{M\otimes H}$ and $\varepsilon _{(M\otimes H)\otimes H}$
are bijective, hence $\varepsilon _{M}$ is too, by the Five Lemma.

We move now to the unit $u_{(-)}$ of the adjunction. Let $N$ be a right $B$%
-module. Consider the maps $i_{1},i_{2}:N\otimes _{B}A\longrightarrow
N\otimes _{B}A\otimes _{B}A$, $i_{1}(n\otimes _{B}a)=n\otimes
_{B}1_{A}\otimes _{B}a$, $i_{2}(n\otimes _{B}a)=n\otimes _{B}a\otimes
_{B}1_{A}$ and the short sequence 
\begin{equation*}
0\longrightarrow N\longrightarrow N\otimes _{B}A\underset{i_{2}}{\overset{%
i_{1}}{\rightrightarrows }}N\otimes _{B}A\otimes _{B}A
\end{equation*}%
where the first morphism is sending $n$ to $n\otimes _{B}1_{A}$. As $A$ is $%
B $-flat, this map is injective. Although the associativity of $A$ fails,
the faithfully flatness property and the existence of the multiplication and
of the unit for $A$ allow us to show, as in the classical case, the
exactness of the sequence in the middle term $N\otimes _{B}A$ (we tensor
over $B$ one more time with $A$, this is easy to see that is exact, and by
faithfully flatness of $A$ we go back to our sequence). Therefore we may
consider the diagram with the top row exact%
\begin{equation*}
\begin{array}{cccccc}
0\longrightarrow & N & \longrightarrow & N\otimes _{B}A & \underset{i_{2}}{%
\overset{i_{1}}{\rightrightarrows }} & N\otimes _{B}A\otimes _{B}A \\ 
& u_{N}\downarrow &  & \parallel &  & I_{N}\otimes _{B}can\downarrow \\ 
0\longrightarrow & (N\otimes _{B}A)^{coH} & \hookrightarrow & N\otimes _{B}A
& \underset{I_{N}\otimes _{B}I_{A}\otimes u_{H}}{\overset{I_{N}\otimes _{B}%
\overline{\rho }_{A}}{\rightrightarrows }} & N\otimes _{B}A\otimes H%
\end{array}%
\end{equation*}%
\newline

In the bottom row, the map $\overline{\rho }_{A}$ is given by $\overline{%
\rho }_{A}(a)=a_{0}\otimes \beta (a_{1})a_{2}$. We need to check the
exactness of this row. Consider $\sum\limits_{i}n_{i}\otimes _{B}a_{i}\in
(N\otimes _{B}A)^{coH}$. Then $\sum\limits_{i}n_{i}\otimes _{B}a_{i0}\otimes
a_{i1}=\sum\limits_{i}n_{i}\otimes _{B}a_{i}\otimes 1_{H}$ implies%
\begin{eqnarray*}
(I_{N}\otimes _{B}\overline{\rho }_{A}-I_{N}\otimes _{B}I_{A}\otimes
u_{H})(\sum\limits_{i}n_{i}\otimes _{B}a_{i}) &=&\sum\limits_{i}n_{i}\otimes
_{B}a_{i0}\otimes \beta (a_{i1})a_{i2}-\sum\limits_{i}n_{i}\otimes
_{B}a_{i}\otimes 1_{H} \\
&=&0
\end{eqnarray*}%
Conversely, let $\sum\limits_{i}n_{i}\otimes _{B}a_{i}\in N\otimes _{B}A$
which satisfies $\sum\limits_{i}n_{i}\otimes _{B}a_{i0}\otimes \beta
(a_{i1})a_{i2}=\sum\limits_{i}n_{i}\otimes _{B}a_{i}\otimes 1_{H}$, and
apply $\rho _{A}$ and $(\Delta \otimes I_{H}\otimes I_{H})(\Delta \otimes
I_{H})\Delta $ on the second, respectively last component of the tensor
product. We obtain 
\begin{equation*}
\sum\limits_{i}n_{i}\otimes _{B}a_{i0}\otimes a_{i1}\otimes \beta
(a_{i2})a_{i3}\otimes a_{i4}\otimes a_{i5}\otimes
a_{i6}=\sum\limits_{i}n_{i}\otimes _{B}a_{i0}\otimes a_{i1}\otimes
1_{H}\otimes 1_{H}\otimes 1_{H}\otimes 1_{H}
\end{equation*}%
Now act by $S$ and $\alpha $ on the forth, respectively fifth tensorand and
apply $\omega $ on the result. It follows that 
\begin{equation*}
\sum\limits_{i}n_{i}\otimes _{B}a_{i0}\otimes \omega (a_{i1},\beta
(a_{i2})S(a_{i3})\alpha (a_{i4}),a_{i5})a_{i6}=\sum\limits_{i}n_{i}\otimes
_{B}a_{i0}\otimes \omega (a_{i1},S(1_{H})\alpha (1_{H}),1_{H})1_{H}
\end{equation*}%
meaning 
\begin{equation*}
\sum\limits_{i}n_{i}\otimes _{B}a_{i0}\otimes
a_{i1}=\sum\limits_{i}n_{i}\otimes _{B}a_{i}\otimes 1_{H}
\end{equation*}

Therefore, the bottom row is exact. The top row is exact by the previous
remarks, while the commutativity of the whole diagram can be easily checked.
Therefore, $u_{N}$ is bijective by the Five Lemma.

(2) $\Longrightarrow $ (1) Follows as in the Hopf case, using also Corollary
(\ref{counit bij implies Galois map bij}).
\end{proof}

\begin{remark}
\rm%
In the proof of the bijectivity of the counit, we have replaced the Galois
maps with $\varepsilon _{M\otimes H}$ and $\varepsilon _{(M\otimes H)\otimes
H}$. Although in the Hopf algebra case they coincide, in our context the
presence of the twist $\widetilde{\mathbf{f}}$ made very difficult to check
directly the commutativity of both diagrams. Therefore we have chosen a
functorial approach, with appropriately changed morphisms. For the proof of
the bijectivity of the second adjunction map $u_{N}$, where $N\in \mathcal{M}%
_{B}$, a change of morphisms in the horizontal rows was also necessary.
\end{remark}

We are going to prove now an affineness condition for coquasi-Hopf algebras.
First we need the following

\begin{proposition}
Let $H$ be a coquasi-Hopf algebra with bijective antipode, $A$ an $H$%
-comodule algebra and $B=A^{coH}$. Assume that there exists $\gamma
:H\longrightarrow A$ a total integral (i.e. a colinear map satisfying $%
\gamma (1_{H})=1_{A}$). Then $u_{N}:N\longrightarrow (N\otimes _{B}A)^{coH}$%
, $u_{N}(n)=n\otimes _{B}1_{A}$, is an isomorphism of right $B$-modules for
all $N\in \mathcal{M}_{B}$.
\end{proposition}

\begin{proof}
We shall define first an analogue of the trace map, namely $%
t_{A}:A\longrightarrow B$, $t_{A}(a)=a_{0}\beta (a_{1})\gamma S(a_{2})$.
This is well defined, because 
\begin{eqnarray*}
\rho t_{A}(a) &=&a_{0_{0}}\beta (a_{1})\gamma S(a_{2})_{0}\otimes
a_{0_{1}}\gamma S(a_{2})_{1} \\
(\text{colinearity of }\gamma ) &=&a_{0}\beta (a_{2})\gamma S(a_{4})\otimes
a_{1}S(a_{3}) \\
&=&a_{0}\beta (a_{1})\gamma S(a_{2})\otimes 1_{H}
\end{eqnarray*}%
where we have used that $S$ is an antimorphism of coalgebras and relation (%
\ref{IdbetaS}). Then using again relation (\ref{IdbetaS}) one can check that
the map $(N\otimes _{B}A)^{coH}\longrightarrow N$, $\sum_{i}n_{i}\otimes
_{B}a_{i}\longrightarrow n_{i}t_{A}(a_{i})$ is the inverse of $u_{N}$.
\end{proof}

We may state now the announced affineness criterion:

\begin{theorem}
\label{total integral implica strong str th}Let $H$ be a coquasi-Hopf
algebra with bijective antipode, $A$ an $H$-comodule algebra and $B=A^{coH}$%
. Assume that

(1) There exists $\gamma :H\longrightarrow A$ a total integral;

(2) The canonical map $can:A\otimes _{B}A\longrightarrow A\otimes H$ is
surjective.

Then the functor of coinvariants and the induction functor form a pair of
inverse equivalences between $\mathcal{M}_{A}^{H}$ and $\mathcal{M}_{B}$.
\end{theorem}

\begin{proof}
From the previous Proposition, we know that the unit of the adjunction is
bijective. It remains to show that $\varepsilon _{M}:M^{coH}\otimes
_{B}A\longrightarrow M$ is an isomorphism for any Hopf module $M\in \mathcal{%
M}_{A}^{H}$. We shall follow here the approach from \cite{Schauenburg04}.

Recall that Bulacu and Nauwelaerts (\cite{Bulacu00})\ have proven the
equivalence between the existence of a total integral on a comodule algebra $%
A$ and the injectivity of any Hopf module as a right $H$-comodule. Their
result is stated for right Hopf modules, but it holds also for $_{A}\mathcal{%
M}^{H}$ because the antipode is bijective and $\mathcal{M}_{A}^{H}\simeq
{}_{A_{\mathbf{f}^{-1}}^{op}}\mathcal{M}^{H^{cop}}$.

From Remark \ref{can morf de hopf mod} we know that $can$ is a morphism of
left Hopf modules. The composition 
\begin{equation*}
\widetilde{can}:{}_{\bullet }A^{\bullet }\otimes A\longrightarrow
{}_{\bullet }A^{\bullet }\otimes _{B}A\overset{can}{\longrightarrow }%
{}_{\bullet }A^{\bullet }\otimes H^{S,\bullet }
\end{equation*}%
will be a surjective left Hopf module map, therefore it splits as an $H$%
-comodule map via a colinear morphism $\theta :_{\bullet }A^{\bullet
}\otimes H^{S,\bullet }\longrightarrow {}_{\bullet }A^{\bullet }\otimes A$
with $\widetilde{can}\theta =I_{A\otimes H}$. Denote $\theta (1_{A}\otimes
h)=\sum\limits_{i}\widetilde{l}_{i}(h)\otimes \widetilde{r}_{i}(h)$ as an
extension of the notation from Proposition \ref{prop can a lin}. It follows
that analogues of properties (5.1)-(5.4) hold:

\begin{eqnarray}
\sum\limits_{i}\widetilde{l}_{i}(h_{1})\otimes \widetilde{r}%
_{i}(h_{1})\otimes h_{2} &=&\sum\limits_{i}\widetilde{l}_{i}(h)\otimes 
\widetilde{r}_{i}(h)_{0}\otimes \widetilde{r}_{i}(h)_{1}  \label{r0} \\
\sum\limits_{i}\widetilde{l}_{i}(h)\widetilde{r}_{i}(h) &=&\alpha (h)1_{A}
\label{lr} \\
\sum\limits_{i}\widetilde{l}_{i}(h)_{0}\otimes \widetilde{r}_{i}(h)\otimes 
\widetilde{l}_{i}(h)_{1} &=&\sum\limits_{i}\widetilde{l}_{i}(h_{2})\otimes 
\widetilde{r}_{i}(h_{2})\otimes S(h_{1})  \label{l0} \\
\sum\limits_{i}a_{0}\beta (a_{1})\widetilde{l}_{i}(a_{2})\otimes \widetilde{r%
}_{i}(a_{2}) &=&1_{A}\otimes a  \label{alr}
\end{eqnarray}%
Relation $(\ref{l0})$ implies $\sum\limits_{i}m_{0}\beta (m_{1})\widetilde{l}%
_{i}(m_{2})\otimes \widetilde{r}_{i}(m_{2})\in M^{coH}\otimes A$. Now we can
define $\chi _{M}:M\longrightarrow M^{coH}\otimes _{B}A$, $\chi
_{M}(m)=\sum\limits_{i}m_{0}\beta (m_{1})\widetilde{l}_{i}(m_{2})\otimes _{B}%
\widetilde{r}_{i}(m_{2})$. We claim that this is an inverse for $\varepsilon
_{M}$, for any $M\in \mathcal{M}_{A}^{H}$. Indeed%
\begin{eqnarray*}
\chi _{M}\varepsilon _{M}(m\otimes _{B}a) &=&\chi
_{M}(ma)=\sum_{i}ma_{0}\beta (a_{1})\widetilde{l}_{i}(a_{2})\otimes _{B}%
\widetilde{r}_{i}(a_{2}) \\
(\ref{alr}) &=&m\otimes _{B}a
\end{eqnarray*}%
for all $m\otimes _{B}a\in M^{coH}\otimes _{B}A$ and%
\begin{eqnarray*}
\varepsilon _{M}\chi _{M}(m) &=&\varepsilon _{M}(\sum_{i}m_{0}\beta (m_{1})%
\widetilde{l}_{i}(m_{2})\otimes _{B}\widetilde{r}_{i}(m_{2})) \\
&=&\sum_{i}[m_{0}\beta (m_{1})\widetilde{l}_{i}(m_{2})]\widetilde{r}%
_{i}(m_{2}) \\
&=&\sum_{i}m_{0}[\widetilde{l}_{i}(m_{3})_{0}\widetilde{r}%
_{i}(m_{3})_{0}]\omega (m_{1},\widetilde{l}_{i}(m_{3})_{1},\widetilde{r}%
_{i}(m_{3})_{1})\beta (m_{2}) \\
(\ref{l0}) &=&\sum_{i}m_{0}[\widetilde{l}_{i}(m_{4})\widetilde{r}%
_{i}(m_{4})_{0}]\omega (m_{1},S(m_{3}),\widetilde{r}_{i}(m_{4})_{1})\beta
(m_{2}) \\
(\ref{r0}) &=&\sum_{i}m_{0}[\widetilde{l}_{i}(m_{4})\widetilde{r}%
_{i}(m_{4})]\omega (m_{1},S(m_{3}),m_{5})\beta (m_{2}) \\
(\ref{lr}) &=&\sum_{i}m_{0}\omega (m_{1},S(m_{3}),m_{5})\beta (m_{2})\alpha
(m_{4}) \\
(\ref{omega anihileaza S}) &=&m
\end{eqnarray*}%
for all $m\in M$. It follows that $\varepsilon _{M}$ is bijective.
\end{proof}

We can state now all our previous results in the form of the following
theorem:

\begin{theorem}
\label{big}Let $H$ be a coquasi-Hopf algebra with bijective antipode, $A$ an 
$H$-comodule algebra and $B=A^{coH}$. Then the following are equivalent:

(1) There exists a total integral $\gamma :H\longrightarrow A$ and the map $%
can:A\otimes _{B}A\longrightarrow {}A\otimes H$ is surjective;

(2) The functor of coinvariants and the induction functor are a pair of
inverse equivalences between $\mathcal{M}_{A}^{H}$ and $\mathcal{M}_{B}$;

(3) The functor of coinvariants and the induction functor are a pair of
inverse equivalences between $_{A}\mathcal{M}^{H}$ and $_{B}\mathcal{M}$;

(4) $A$ is faithfully flat as a left $B$-module, and $A$ is a Galois
extension of $B$;

(5) $A$ is faithfully flat as a right $B$-module, and $A$ is a Galois
extension of $B$.
\end{theorem}

\begin{proof}
(1) $\Longrightarrow $ (2) follows from Theorem \ref{total integral implica
strong str th}. (2) $\Longleftrightarrow $ (4) is Theorem \ref{a galois +
fidel plat echiv strong str th}. (4) $\Longrightarrow $ (1) uses the same
argument as in \cite{Schneider90a}, because $can^{\prime }$ is also
bijective by Lemma \ref{can bij implica canprim bij}, and is a morphism of
left $B$-modules, right $H$-comodules by Proposition \ref{schneider can-1}%
.(6). The sequence of isomorphisms ($A$ is flat $B$-module) 
\begin{equation*}
(A^{\bullet }{\small \square }_{H}V)\otimes _{B}A\simeq (A^{\bullet }\otimes
_{B}A){\small \square }_{H}V\simeq (A\otimes H^{\bullet }){\small \square }%
_{H}V\simeq A\otimes V
\end{equation*}%
for each $V\in {}^{H}\mathcal{M}$, together with the left $B$-faithful
flatness of $A$ imply that $A$ is right $H$-coflat, or equivalently, that $A$
is $H$-injective (here ${\small \square }_{H}$ is the cotensor product over $%
H$).

(1) $\Longleftrightarrow $ (3) $\Longleftrightarrow $ (5) We simply apply
the above to $A^{op}$ as a right $H^{op}$-comodule algebra, since the
antipode is bijective.
\end{proof}

\section{A bialgebroid associated to a faithfully flat Galois extension}

Let $H$ be a coquasi-bialgebra (without any assumption on the antipode) and $%
A$ a right $H$-comodule algebra. On the tensor product $A\otimes A^{op}$ we
consider the codiagonal right $H$-comodule structure $\rho (a\otimes
b)=a_{0}\otimes b_{0}\otimes a_{1}b_{1}$. Denote $L=(A\otimes A^{op})^{coH}$%
. Then

\begin{proposition}
$L$ is an associative $B\otimes B^{op}$-algebra with unit $1_{A}\otimes
1_{A} $ and multiplication 
\begin{equation}
(a\otimes b)(c\otimes d)=a_{0}c_{0}\otimes d_{0}b_{0}\omega
^{-1}(a_{1},c_{1},d_{1}b_{1})\omega (c_{2},d_{2},b_{2})
\label{multilplicarea pe (AtensorA)coH}
\end{equation}%
for $a\otimes b$, $c\otimes d\in L$.
\end{proposition}

\begin{proof}
We shall suppress the $\sum $ symbol when referring to elements of $L$ for
simplicity.

After a short calculation, it follows that the multiplication is
well-defined, with values in $L$. Moreover, the maps $b\in B\longrightarrow
b\otimes 1_{A}\in A\otimes A^{op}$, $b\in B^{op}\longrightarrow 1_{A}\otimes
b\in A\otimes A^{op}$ take values in $L$ and are multiplicative.

It is easy to check that $1_{A}\otimes 1_{A}\in L$ and that it is a unit for
the given multiplication. The most difficult part to show is the
associativity. Take $a\otimes b$, $c\otimes d$, $e\otimes f\in L$ (summation
understood). Then we compute%
\begin{allowdisplaybreaks}
\begin{eqnarray*}
\lbrack (a\otimes b)(c\otimes d)](e\otimes f) &=&(a_{0}c_{0}\otimes
d_{0}b_{0})(e\otimes f)\omega ^{-1}(a_{1},c_{1},d_{1}b_{1})\omega
(c_{2},d_{2},b_{2}) \\
(\ref{multilplicarea pe (AtensorA)coH}) &=&(a_{0}c_{0})e_{0}\otimes
f_{0}(d_{0}b_{0})\omega ^{-1}(a_{1}c_{1},e_{1},f_{1}(d_{1}b_{1}))\omega
(e_{2},f_{2},d_{2}b_{2}) \\
&&\omega ^{-1}(a_{2},c_{2},d_{3}b_{3})\omega (c_{3},d_{4},b_{4}) \\
(\ref{asociat multipl}) &=&(a_{0}c_{0})e_{0}\otimes (f_{0}d_{0})b_{0}\omega
^{-1}(a_{1}c_{1},e_{1},(f_{1}d_{1})b_{1})\omega ^{-1}(f_{2},d_{2},b_{2}) \\
&&\omega (e_{2},f_{3},d_{3}b_{3})\omega ^{-1}(a_{2},c_{2},d_{4}b_{4})\omega
(c_{3},d_{5},b_{5}) \\
(\ref{cocycle omega}) &=&(a_{0}c_{0})e_{0}\otimes (f_{0}d_{0})b_{0}\omega
^{-1}(a_{1}c_{1},e_{1},(f_{1}d_{1})b_{1})\omega (e_{2},f_{2}d_{2},b_{2}) \\
&&\omega (e_{3},f_{3},d_{3})\omega ^{-1}(e_{4}f_{4},d_{4},b_{3})\omega
^{-1}(a_{2},c_{2},d_{5}b_{4})\omega (c_{3},d_{6},b_{5}) \\
(e\otimes f\in L) &=&(a_{0}c_{0})e_{0}\otimes (f_{0}d_{0})b_{0}\omega
^{-1}(a_{1}c_{1},e_{1},(f_{1}d_{1})b_{1})\omega (e_{2},f_{2}d_{2},b_{2}) \\
&&\omega (e_{3},f_{3},d_{3})\omega ^{-1}(a_{2},c_{2},d_{4}b_{3})\omega
(c_{3},d_{5},b_{4}) \\
(\ref{cocycle omega}) &=&(a_{0}c_{0})e_{0}\otimes (f_{0}d_{0})b_{0}\omega
((a_{1}c_{1})e_{1},f_{1}d_{1},b_{1})\omega ^{-1}(a_{2}c_{2},e_{2},f_{2}d_{2})
\\
&&\omega ^{-1}(a_{3}c_{3},e_{3}(f_{3}d_{3}),b_{2})\omega
(e_{4},f_{4},d_{4})\omega ^{-1}(a_{4},c_{4},d_{5}b_{3}) \\
&&\omega (c_{5},d_{6},b_{4}) \\
(\ref{asociat multipl}),(e\otimes f\in L) &=&(a_{0}c_{0})e_{0}\otimes
(f_{0}d_{0})b_{0}\omega ((a_{1}c_{1})e_{1},f_{1}d_{1},b_{1})\omega
^{-1}(a_{2}c_{2},e_{2},f_{2}d_{2}) \\
&&\omega (e_{3},f_{3},d_{3})\omega ^{-1}(a_{3}c_{3},d_{4},b_{2})\omega
^{-1}(a_{4},c_{4},d_{5}b_{3})\omega (c_{5},d_{6},b_{4}) \\
(\ref{cocycle omega}),(c\otimes d\in L) &=&(a_{0}c_{0})e_{0}\otimes
(f_{0}d_{0})b_{0}\omega ((a_{1}c_{1})e_{1},f_{1}d_{1},b_{1})\omega
^{-1}(a_{2}c_{2},e_{2},f_{2}d_{2}) \\
&&\omega (e_{3},f_{3},d_{3})\omega ^{-1}(a_{3},c_{3},d_{4}) \\
(\ref{asociat multipl}) &=&a_{0}(c_{0}e_{0})\otimes (f_{0}d_{0})b_{0}\omega
(a_{1}(c_{1}e_{1}),f_{1}d_{1},b_{1})\omega (a_{2},c_{2},e_{2}) \\
&&\omega ^{-1}(a_{3}c_{3},e_{3},f_{2}d_{2})\omega (e_{4},f_{3},d_{3})\omega
^{-1}(a_{4},c_{4},d_{4})
\end{eqnarray*}%
\end{allowdisplaybreaks}%

Apparently this leads nowhere. But let's evaluate also%
\begin{allowdisplaybreaks}%
\begin{eqnarray*}
(a\otimes b)[(c\otimes d)(e\otimes f)] &=&(a\otimes b)(c_{0}e_{0}\otimes
f_{0}d_{0})\omega ^{-1}(c_{1},e_{1},f_{1}d_{1})\omega (e_{2},f_{2},d_{2}) \\
(\ref{multilplicarea pe (AtensorA)coH}) &=&a_{0}(c_{0}e_{0})\otimes
(f_{0}d_{0})b_{0}\omega ^{-1}(a_{1},c_{1}e_{1},(f_{1}d_{1})b_{1})\omega
(c_{2}e_{2},f_{2}d_{2},b_{2}) \\
&&\omega ^{-1}(c_{3},e_{3},f_{3}d_{3})\omega (e_{4},f_{4},d_{4}) \\
(\ref{cocycle omega}) &=&a_{0}(c_{0}e_{0})\otimes (f_{0}d_{0})b_{0}\omega
(a_{1}(c_{1}e_{1}),f_{1}d_{1},b_{1})\omega ^{-1}(a_{2},c_{2}e_{2},f_{2}d_{2})
\\
&&\omega ^{-1}(a_{3},(c_{3}e_{3})(f_{3}d_{3}),b_{2})\omega
^{-1}(c_{4},e_{4},f_{4}d_{4})\omega (e_{5},f_{5},d_{5}) \\
(\ref{asociat multipl}),(\ref{asociat multipl}) &=&a_{0}(c_{0}e_{0})\otimes
(f_{0}d_{0})b_{0}\omega (a_{1}(c_{1}e_{1}),f_{1}d_{1},b_{1})\omega
^{-1}(a_{2},c_{2}e_{2},f_{2}d_{2}) \\
&&\omega ^{-1}(c_{3},e_{3},f_{3}d_{3})\omega (e_{4},f_{4},d_{4})\omega
^{-1}(a_{3},c_{4}d_{5},b_{2}) \\
(c\otimes d\in L) &=&a_{0}(c_{0}e_{0})\otimes (f_{0}d_{0})b_{0}\omega
(a_{1}(c_{1}e_{1}),f_{1}d_{1},b_{1})\omega ^{-1}(a_{2},c_{2}e_{2},f_{2}d_{2})
\\
&&\omega ^{-1}(c_{3},e_{3},f_{3}d_{3})\omega (e_{4},f_{4},d_{4}) \\
(\ref{cocycle omega}) &=&a_{0}(c_{0}e_{0})\otimes (f_{0}d_{0})b_{0}\omega
(a_{1}(c_{1}e_{1}),f_{1}d_{1},b_{1})\omega (a_{2},c_{2},e_{2}) \\
&&\omega ^{-1}(a_{3}c_{3},e_{3},f_{2}d_{2})\omega
^{-1}(a_{4},c_{4},e_{4}(f_{3}d_{3}))\omega (e_{5},f_{4},d_{4}) \\
(\ref{asociat multipl}), &=&a_{0}(c_{0}e_{0})\otimes (f_{0}d_{0})b_{0}\omega
(a_{1}(c_{1}e_{1}),f_{1}d_{1},b_{1})\omega (a_{2},c_{2},e_{2}) \\
&&\omega ^{-1}(a_{3}c_{3},e_{3},f_{2}d_{2})\omega (e_{4},f_{3},d_{3})\omega
^{-1}(a_{4},c_{4},d_{4})
\end{eqnarray*}%
\end{allowdisplaybreaks}%
hence we have obtained the same as above.
\end{proof}

\begin{proposition}
\label{prop bialgebroid}Let $H$ be a coquasi-Hopf algebra with bijective
antipode and $A$ a right $H$-comodule algebra, left faithfully flat and
Galois over $B=A^{coH}$. Then the left $L$-module category $_{L}\mathcal{M}$
is equivalent to the category of two sided $(H,A)$-Hopf modules $_{A}%
\mathcal{M}_{A}^{H}$(i.e. $A$-bimodules in $\mathcal{M}^{H}$).
\end{proposition}

\begin{proof}
Let $N$ be any left $L$-module. By restriction, $N$ is a left $B^{op}$%
-module, that is, a right $B$-module. We can use then the category
equivalence $\mathcal{M}_{B}\simeq \mathcal{M}_{A}^{H}$ from Theorem \ref{a
galois + fidel plat echiv strong str th}. It follows that $N\otimes
_{B}A_{\bullet }^{\bullet }\in \mathcal{M}_{A}^{H}$. But we still need the
left $A$-module structure on $N\otimes _{B}A_{\bullet }^{\bullet }$. For
this, we shall use the inverse of the Galois map with notations from
Proposition \ref{prop can a lin}:%
\begin{equation}
a\odot (n\otimes _{B}b)=\sum_{i}[a_{0}\otimes b_{0}\omega (a_{1},b_{1},\beta
(a_{2}b_{2})S(a_{3}b_{3}))l_{i}(a_{4}b_{4})]n\otimes _{B}r_{i}(a_{4}b_{4})
\label{a modul stang pe indus}
\end{equation}%
for any $a,b\in A$, $n\in N$. We should check first if this is well-defined.
For this, consider $(A^{\bullet }\otimes A^{\bullet })\otimes _{B}A$ as a
right $H$-comodule with coaction on the first component. We can compute%
\begin{eqnarray*}
\rho ([a_{0}\otimes b_{0}\omega (a_{1},b_{1},\beta
(a_{2}b_{2})S(a_{3}b_{3}))l_{i}(a_{4}b_{4})]\otimes _{B}r_{i}(a_{4}b_{4}))
&=&\sum_{i}[a_{0}\otimes b_{0}\omega (a_{2},b_{2},\beta
(a_{3}b_{3})S(a_{4}b_{4})) \\
&&l_{i}(a_{5}b_{5})_{0}]\otimes _{B}r_{i}(a_{6}b_{6})\otimes
a_{1}[b_{1}l_{i}(a_{5}b_{5})_{1}] \\
(\text{Proposition }\ref{prop can a lin}(5.3)),(\ref{asociat multipl}),(\ref%
{IdbetaS}) &=&\sum_{i}[a_{0}\otimes b_{0}\omega (a_{1},b_{1},\beta
(a_{2}b_{2})S(a_{3}b_{3})) \\
&&l_{i}(a_{4}b_{4})]\otimes _{B}r_{i}(a_{4}b_{4})\otimes 1_{H}
\end{eqnarray*}%
By the left $B$-flatness of $A$, it follows that 
\begin{equation*}
\sum\limits_{i}[a_{0}\otimes b_{0}\omega (a_{1},b_{1},\beta
(a_{2}b_{2})S(a_{3}b_{3}))l_{i}(a_{4}b_{4})]\otimes _{B}r_{i}(a_{4}b_{4})\in
(A^{\bullet }\otimes A^{\bullet })^{coH}\otimes _{B}A
\end{equation*}%
Therefore, relation (\ref{a modul stang pe indus}) is correct. We still need
to check that the formula defines on $N\otimes _{B}A_{\bullet }^{\bullet }$
an $A$-bimodule structure in $\mathcal{M}^{H}$. The left $A$-multiplication
is $H$-colinear and 
\begin{eqnarray*}
1_{A}\odot (n\otimes _{B}b) &=&\sum_{i}[1_{A}\otimes b_{0}\beta
(b_{1})l_{i}(b_{2})]n\otimes _{B}r_{i}(b_{2}) \\
(\text{Proposition \ref{prop can}(5.4)}) &=&(1_{A}\otimes 1_{A})n\otimes
_{B}b \\
&=&n\otimes _{B}b
\end{eqnarray*}%
Now we compute%
\begin{allowdisplaybreaks}
\begin{eqnarray*}
a_{0}\odot (b_{0}\odot (n\otimes _{B}c_{0}))\omega (a_{1},b_{1},c_{1})
&=&\sum_{i}a_{0}\odot \lbrack \lbrack b_{0}\otimes c_{0}\omega
(b_{1},c_{1},\beta (b_{2}c_{2})S(b_{3}c_{3})) \\
&&l_{i}(b_{4}c_{4})]n\otimes _{B}r_{i}(b_{4}c_{4})]\omega (a_{1},b_{5},c_{5})
\\
(\text{Proposition \ref{prop can}(5.1)}) &=&\sum_{i,j}(a_{0}\otimes
r_{i}(b_{4}c_{4})l_{j}(a_{4}(b_{8}c_{8})))(b_{0}\otimes
c_{0}l_{i}(b_{4}c_{4}))n \\
&&\otimes _{B}r_{j}(a_{4}(b_{8}c_{8}))\omega
(a_{1},b_{5}c_{5},S(a_{3}(b_{7}c_{7})) \\
&&\omega (b_{1},c_{1},S(b_{3}c_{3}))\beta (b_{2}c_{2})\beta
(a_{2}(b_{6}c_{6})) \\
(\ref{multilplicarea pe (AtensorA)coH}),(\text{Proposition \ref{prop can}%
(5.1),(5.3)}) &=&\sum_{i,j}\left\{ a_{0}b_{0}\otimes \lbrack
c_{0}l_{i}(b_{8}c_{8})][r_{i}(b_{8}c_{8})l_{j}(a_{7}(b_{16}c_{16}))]\right\}
n \\
&&\otimes _{B}r_{j}(a_{7}(b_{16}c_{16}))\omega
^{-1}(a_{1},b_{1},[c_{1}S(b_{7}c_{7})][(b_{9}c_{9}) \\
&&S(a_{6}(b_{15}c_{15}))])\omega (b_{2},c_{2}S(b_{6}c_{6}),(b_{10}c_{10}) \\
&&S(a_{5}(b_{14}c_{14})))\omega (a_{2},b_{11}c_{11},S(a_{4}(b_{13}c_{13})))
\\
&&\omega (b_{3},c_{3},S(b_{5}c_{5}))\beta (b_{4}c_{4})\beta
(a_{3}(b_{12}c_{12})) \\
(\ref{cocycle omega}),(\ref{asociat multipl}),(\ref{asoc comod alg})
&=&\sum_{i,j}\left\{ a_{0}b_{0}\otimes
c_{0}[l_{i}(b_{8}c_{8})[r_{i}(b_{8}c_{8})l_{j}(a_{8}(b_{18}c_{18}))]]\right%
\} n \\
&&\otimes _{B}r_{j}(a_{8}(b_{18}c_{18}))\omega
^{-1}(a_{1},b_{1},c_{1}[S(b_{7}c_{7})[(b_{9}c_{9}) \\
&&S(a_{7}(b_{17}c_{17}))]])\omega (b_{2},c_{2},S(b_{6}c_{6})[(b_{10}c_{10})
\\
&&S(a_{6}(b_{16}c_{16}))])\omega (b_{3}c_{3},S(b_{5}c_{5}),(b_{11}c_{11}) \\
&&S(a_{5}(b_{15}c_{15})))\omega (a_{2},b_{12}c_{12},S(a_{4}(b_{14}c_{14})))
\\
&&\beta (b_{4}c_{4})\beta (a_{3}(b_{13}c_{13})) \\
(\ref{asociat multipl}),(\ref{asoc comod alg}),(\text{Proposition \ref{prop
can}(5.2)}),(\ref{SalfaId}) &=&\sum_{j}\left\{ a_{0}b_{0}\otimes
c_{0}l_{j}(a_{9}(b_{17}c_{17}))\right\} n\otimes
_{B}r_{j}(a_{9}(b_{17}c_{17})) \\
&&\omega ^{-1}(a_{1},b_{1},c_{1}S(a_{8}(b_{16}c_{16}))) \\
&&\omega (b_{2},c_{2},S(a_{7}(b_{15}c_{15}))) \\
&&\omega ^{-1}(S(b_{6}c_{6}),b_{8}c_{8},S(a_{6}(b_{14}c_{14})))\omega
(b_{3}c_{3},S(b_{5}c_{5}), \\
&&(b_{9}c_{9})S(a_{5}(b_{13}c_{13})))\omega
(a_{2},b_{10}c_{10},S(a_{4}(b_{12}c_{12}))) \\
&&\beta (b_{4}c_{4})\beta (a_{3}(b_{11}c_{11})) \\
(\ref{cocycle omega}),(\ref{SalfaId}),(\ref{IdbetaS}),(\ref{omega anihileaza
S}) &=&\sum_{j}\left\{ a_{0}b_{0}\otimes
c_{0}l_{j}(a_{7}(b_{8}c_{8}))\right\} n\otimes _{B}r_{j}(a_{7}(b_{8}c_{8}))
\\
&&\omega ^{-1}(a_{1},b_{1},c_{1}S(a_{6}(b_{7}c_{7})))\omega
(b_{2},c_{2},S(a_{5}(b_{6}c_{6}))) \\
&&\omega (a_{2},b_{3}c_{3},S(a_{4}(b_{5}c_{5})))\beta (a_{3}(b_{4}c_{4})) \\
(\ref{cocycle omega}),(\ref{asociat multipl}) &=&\sum_{j}\left\{
a_{0}b_{0}\otimes c_{0}l_{j}((a_{4}b_{4})c_{4})\right\} n\otimes
_{B}r_{j}((a_{4}b_{4})c_{8}) \\
&&\omega (a_{1}b_{1},c_{1}S((a_{3}b_{3})c_{3})\beta ((a_{2}b_{2})c_{2}) \\
&=&(ab)\odot (n\otimes _{B}c)
\end{eqnarray*}%
\end{allowdisplaybreaks}%
hence $N\otimes _{B}A$ is a left $A$-module in $\mathcal{M}^{H}$. We only
have to show the compatibility between the two $A$-module structures:%
\begin{allowdisplaybreaks}%
\begin{eqnarray*}
a_{0}\odot \lbrack (n\otimes _{B}b_{0})c_{0}]\omega (a_{1},b_{1},c_{1})
&=&\sum_{i}[a_{0}\otimes (b_{0}c_{0})l_{i}(a_{4}(b_{4}c_{4}))]n\otimes
_{B}r_{i}(a_{4} \\
&&(b_{4}c_{4}))\omega (a_{1},b_{1}c_{1},S(a_{3}(b_{3}c_{3}))\omega
(a_{5},b_{5},c_{5}) \\
&&\beta (a_{2}(b_{2}c_{2})) \\
(\ref{asociat multipl}),(\text{Proposition \ref{prop can}(5.3),(5.5)}),(\ref%
{twist f}),(\ref{beta*f-1=delta}) &=&\sum_{i,j}[a_{0}\otimes
(b_{0}c_{0})l_{i}(c_{5})l_{j}(a_{5}b_{5})]n\otimes _{B}r_{j}(a_{5}b_{5}) \\
&&r_{i}(c_{5})\omega (a_{1},b_{1}c_{1},S(c_{4})S(a_{4}b_{4}))\omega
(a_{2},b_{2},c_{2}) \\
&&q(a_{3}b_{3},c_{3}) \\
(\ref{formula q}),(\ref{asoc comod alg}),(\ref{asoc comod alg}),(\text{%
Proposition \ref{prop can}(5.3)}) &=&\sum_{i,j}[a_{0}\otimes b_{0}\left\{
[c_{0}l_{i}(c_{13})]l_{j}(a_{10}b_{11})\right\} ]n\otimes _{B} \\
&&r_{j}(a_{10}b_{11})r_{i}(c_{13})\omega
^{-1}(c_{1},S(c_{12}),S(a_{9}b_{10})) \\
&&\omega (b_{1},c_{2},S(c_{11})S(a_{8}b_{9})) \\
&&\omega (a_{1},b_{2}c_{3},S(c_{10})S(a_{7}b_{8})) \\
&&\omega (a_{2},b_{3},c_{4})\omega ((a_{3}b_{4})c_{5},S(c_{9}),S(a_{6}b_{7}))
\\
&&\omega ^{-1}(a_{4}b_{5},c_{6},S(c_{8}))\beta (a_{5}b_{6})\beta (c_{7}) \\
(\ref{cocycle omega}),(\ref{cocycle omega}),(\ref{asociat multipl}),(\ref%
{IdbetaS}) &=&\sum_{i,j}[a_{0}\otimes b_{0}\left\{
[c_{0}l_{i}(c_{2})]l_{j}(a_{4}b_{4})\right\} ]n\otimes _{B} \\
&&r_{j}(a_{4}b_{4})r_{i}(c_{2})\omega (a_{1},b_{1},S(a_{3}b_{3})) \\
&&\beta (a_{2}b_{2})\beta (c_{1}) \\
(\text{Proposition \ref{prop can}(5.4)}) &=&(a\odot (n\otimes _{B}b))c
\end{eqnarray*}%
\end{allowdisplaybreaks}%
Hence $N\otimes _{B}A\in $ $_{A}\mathcal{M}_{A}^{H}$. It is easy to see that
a map $L$-linear $\eta :N_{1}\longrightarrow N_{2}$ induces a morphism$\
\eta \otimes _{B}I_{A}$ in $_{A}\mathcal{M}_{A}^{H}$. We get then a functor $%
\mathcal{F}:{}_{L}\mathcal{M}\longrightarrow {}_{A}\mathcal{M}_{A}^{H}$. For
the inverse construction, let $M\in {}_{A}\mathcal{M}_{A}^{H}$. Then $%
M^{coH}\in \mathcal{M}_{B}={}_{B^{op}}\mathcal{M}$. For any $m\in M^{coH}$
and $a\otimes b\in L$, we may define 
\begin{equation*}
(a\otimes b)m=a(mb)\overset{m\in M^{coH}}{=}(am)b
\end{equation*}%
Using this multiplication, $M^{coH}\in {}_{L}\mathcal{M}$ and we have a
functor $\mathcal{G}:{}_{A}\mathcal{M}_{A}^{H}\longrightarrow {}_{L}\mathcal{%
M}$.

Notice that these two functors are obtained simply restricting the ones in
Theorem \ref{a galois + fidel plat echiv strong str th}. The unit and the
counit are easily checked to be morphisms in the restricted categories.
Therefore we get the category equivalence $_{L}\mathcal{M}\simeq ${}$_{A}%
\mathcal{M}_{A}^{H}$.
\end{proof}

\begin{corollary}
\label{corolar bialgebroid}The category $_{L}\mathcal{M}$ is monoidal.
\end{corollary}

\begin{proof}
As $_{A}\mathcal{M}_{A}^{H}$ is monoidal with $\bigcirc _{A}$ the tensor
product over $A$ in the comodule category, it remains only to transport the
monoidal structure.
\end{proof}

\begin{remark}
\rm%
The previous theorem generalizes Schauenburg's result in the Hopf algebra
case (\cite{Schauenburg98}). In \cite{Schauenburg03}, he gave a categorical
proof, using actions of monoidal categories. All his arguments were purely
categorical, explaining why Schauenburg's construction can also be performed
for coquasi-Hopf algebras. But in order to avoid long and tedious
computations, we preferred the direct approach.

In \cite{Schauenburg98} it was also shown that there is a $\times _{B}$%
-bialgebra (in the sense of Takeuchi) structure on $L$, using precisely the
monoidal structure given by the above corollary. Let see now that a similar
result holds also in the coquasi case. But first we have an inconvenient: we
cannot tensor over $A$, as this is not an associative algebra. This can be
avoided by considering suitable tensor product, namely in the monoidal
category of comodules.
\end{remark}

\begin{lemma}
Let $H$ be a coquasi-bialgebra, $A$ a right $H$-comodule algebra and $%
B=A^{coH}$. For any right $B$-module $N$ and any left Hopf module $M\in
{}_{A}\mathcal{M}^{H}$, we have $(N\otimes _{B}A_{\bullet }^{\bullet
})\bigcirc _{A}M^{\bullet }\simeq N\otimes _{B}M^{\bullet }$ as comodules,
where $\bigcirc _{A}$ denotes the tensor product over $A$ in the monoidal
category $\mathcal{M}^{H}$, $N\otimes _{B}A_{\bullet }^{\bullet }$ is the
induced right Hopf module and $N\otimes _{B}M^{\bullet }$ carries the
comodule structure given by that of $M$.
\end{lemma}

\begin{proof}
Recall that the tensor product over $A$ is the equalizer (in the category of
right comodules) of the following morphisms $j_{1},j_{2}:[(N\otimes
_{B}A)\otimes A]\otimes M\longrightarrow (N\otimes _{B}A)\otimes M$, where%
\begin{eqnarray*}
&&j_{1}([(n\otimes _{B}a)\otimes b]\otimes m=n\otimes _{B}ab\otimes m \\
&&j_{1}([(n\otimes _{B}a)\otimes b]\otimes m=n\otimes _{B}a_{0}\otimes
b_{0}m_{0}\omega (a_{1},b_{1},m_{1})
\end{eqnarray*}%
Now define $\varphi :(N\otimes _{B}A)\otimes M\longrightarrow N\otimes _{B}M$%
, $\varphi ((n\otimes _{B}a)\otimes m)=n\otimes _{B}am$. Then $\varphi $ is
colinear and $\varphi j_{1}=\varphi j_{2}$. Hence it induces the desired
isomorphism.
\end{proof}

\begin{corollary}
\label{corolar monoidal functor}Let $H$ be a coquasi-Hopf algebra with
bijective antipode and $A$ a right $H$-comodule algebra, left faithfully
flat and Galois over $B=A^{coH}$. Then the equivalence $\mathcal{M}%
_{A}^{H}\simeq \mathcal{M}_{B}$ induces a monoidal functor $(-)^{coH}:{}_{A}%
\mathcal{M}_{A}^{H}\simeq {}_{B}\mathcal{M}_{B}$.
\end{corollary}

\begin{proof}
It follows by the previous Lemma and from \cite{Schauenburg98}, Lemma 6.1$.$
\end{proof}

As the monoidal structure of $_{L}\mathcal{M}$ comes from the one of $_{A}%
\mathcal{M}_{A}^{H}$ and the functorial diagram 
\begin{equation*}
\begin{array}{ccc}
_{A}\mathcal{M}_{A}^{H} & \rightleftarrows & _{L}\mathcal{M} \\ 
(-)^{coH}\searrow &  & \swarrow \mathcal{U} \\ 
& _{B^{op}}\mathcal{M} & 
\end{array}%
\end{equation*}%
commutes, where $\mathcal{U}$ is the forgetful functor, from Corollary \ref%
{corolar bialgebroid} it follows that there $\mathcal{U}$ is also monoidal.
But according to \cite{Schauenburg98} and \cite{Brzezinski02a}, a $B\otimes
B^{op}$-algebra $L$ such that the forgetful functor $_{L}\mathcal{%
M\longrightarrow {}}_{B}\mathcal{M}_{B}$ is (strictly) monoidal is precisely
a $\times _{B}$-bialgebra (in the sense of Takeuchi) or equivalently, a
bialgebroid. Therefore we have obtained a new structure object $L$, whose
properties (mainly for the case $B=\Bbbk $) will make the purpose of an
author's forthcoming paper. Having in mind the Hopf algebra case, where the
biGalois extensions and torsors are involved, it is expected that this will
clarify more about the connections between various generalizations of Hopf
algebras.

\begin{acknowledgement}
The author would like to thank Prof. C. N\u{a}st\u{a}sescu and F. Panaite
for their useful comments which improved this paper.
\end{acknowledgement}


\begin{thebibliography}{99}
\bibitem{Albuquerque99} H.~Albuquerque and S.~Majid. \newblock Quasialgebra {%
S}tructure of the {O}ctonions. \newblock {\em J. Algebra}, 220(1):188--224,
1999.

\bibitem{Albuquerque99a} H.~Albuquerque and S.~Majid. \newblock$%
\mathbb{Z}
_{n}$-{Q}uasialgebras. \newblock In \emph{Matrices and Group
Representations. Proceedings of a Workshop Dedicated to Professor Graciano
N. De Oliveira on the Occasion of his 60th Birthday, Coimbra, Portugal, May
6--8, 1998}, volume~19 of \emph{Textos Math., Ser.}, pages 57--64. Coimbra:
Univ. de Coimbra, Departamento de Matem\'{a}tica, 1999.

\bibitem{Altschuler92} D.~Altschuler and A.~Coste. \newblock Quasi-{Q}uantum
groups, {K}nots, {T}hree {M}anifolds and {T}opological {F}ield {T}heory. %
\newblock {\em Comm. Math. Phys.}, 150:83--107, 1992.

\bibitem{Ardizzoni07} A.~Ardizzoni, C.~Menini, and D.~Stefan. \newblock %
Hochschild {C}ohomology and {S}moothness in {M}onoidal {C}ategories. %
\newblock {\em J. Pure Appl. Algebra}, 208:297--330, 2007.

\bibitem{Balan07a} A.~Balan. \newblock Crossed {P}roducts for {C}oquasi-{H}%
opf {A}lgebras. \newblock preprint 2007.

\bibitem{Balan06} A.~Balan. \newblock A {M}orita {C}ontext and {G}alois {E}%
xtensions for {Q}uasi-{H}opf {A}lgebras. \newblock to appear in Comm. Alg.;
available on arXiv:math.QA/0705.3515, 2007.

\bibitem{Blattner89} R.~J. Blattner and S.~Montgomery. \newblock Crossed {P}%
roducts and {G}alois {E}xtensions of {H}opf {A}lgebras. 
\newblock {\em
Pacific J. Math.}, 137:37--54, 1989.

\bibitem{Bohm04} G.~B{\"o}hm. \newblock Galois {T}heory for {H}opf {A}%
lgebroids. \newblock {\em Ann. Univ. Ferrara, Nuova Ser., Sez. VII},
51:233--262, 2005.

\bibitem{Brzezinski99} T.~Brz{\`e}zinski and P.~M. Hajac. \newblock %
Coalgebra {E}xtensions and {A}lgebra {C}oextensions of {G}alois {T}ype. %
\newblock {\em Comm. Algebra}, 27:1347--1367, 1999.

\bibitem{Brzezinski02a} T.~Brzezinski and G.~Militaru. \newblock %
Bialgebroids, $\times_{A}$-{B}ialgebras and {D}uality. 
\newblock {\em J.
Algebra}, 251(1):279--294, 2002.

\bibitem{Bulacu99} D.~Bulacu. \newblock On the {A}ntipode of {S}emi-{H}opf {A%
}lgebras and {B}raided {S}emi-{H}opf {A}lgebras. 
\newblock {\em Rev. Roum.
Math. Pures Appl.}, 44(3):329--340, 1999.

\bibitem{BulacuChirita01} D.~Bulacu and B.~Chirita. \newblock Dual {D}%
rinfeld {D}ouble by {D}iagonal {C}rossed {P}roduct. 
\newblock {\em Rev.
Roum. Math. Pures Appl.}, 47(3):271--294, 2002.

\bibitem{Bulacu00} D.~Bulacu and E.~Nauwelaerts. \newblock Relative {H}opf {M%
}odules for ({D}ual) {Q}uasi-{H}opf {A}lgebras. \newblock {\em J. Algebra},
229(2):632--659, 2000.

\bibitem{Bulacu02co} D.~Bulacu and E.~Nauwelaerts. \newblock Dual {Q}uasi-{H}%
opf {A}lgebra {C}oactions, {S}mash {C}oproducts and {R}elative {H}opf {M}%
odules. \newblock {\em Rev. Roum. Math. Pures Appl.}, 47(4):415--443, 2002.

\bibitem{Chase65} S.~U. Chase, D.~K. Harrison, and A.~Rosenberg. \newblock
\emph{Galois {T}heory and {C}ohomology of {C}ommutative {R}ings}. \newblock %
Number~52 in {AMS} Memoirs. AMS, Providence, R. I., 1965.

\bibitem{Chase69} S.~U. Chase and M.~E. Sweedler. 
\newblock {\em Hopf
{A}lgebras and {G}alois {T}heory}. \newblock Number~97 in Lect. Notes in
Math. Springer Verlag, Berlin, 1969.

\bibitem{Doi86} Y.~Doi and M.~Takeuchi. \newblock Cleft {C}omodule {A}%
lgebras for a {B}ialgebra. \newblock {\em Comm. Alg.}, 14:801--818, 1986.

\bibitem{Doi89} Y.~Doi and M.~Takeuchi. \newblock Hopf-{G}alois {E}xtensions
of {A}lgebras, the {M}iyashita-{U}lbrich {A}ction, and {A}zumaya {A}lgebras. %
\newblock {\em J. Algebra}, 121(2):488--516, 1989.

\bibitem{Drinfeld90} V.~G. Drinfeld. \newblock Quasi-{H}opf {A}lgebras. %
\newblock {\em Leningrad Math. J.}, 1:1419--1457, 1990.

\bibitem{Dascalescu00} S.~D\u{a}sc\u{a}lescu, C.~N\u{a}st\u{a}sescu, and 
\c{S}. Raianu. \newblock {\em Hopf {A}lgebras: {A}n {I}ntroduction}. %
\newblock Number 235 in Pure and Applied Math. Marcel Dekker, New York, 2001.

\bibitem{Kadison03} L.~Kadison and K.~Szlachanyi. \newblock Bialgebroid {A}%
ctions on {D}epth {T}wo {E}xtensions and {D}uality. 
\newblock {\em Adv.
Math.}, 179:75--121, 2003.

\bibitem{Kassel95} C.~Kassel. \newblock {\em Quantum {G}roups}. \newblock %
Number 155 in Graduate Texts in Math. Springer Verlag, 1995.

\bibitem{Kreimer81} H.~F. Kreimer and M.~Takeuchi. \newblock Hopf {A}lgebras
and {G}alois {E}xtensions of an {A}lgebra. \newblock {\em Indiana Math. J.},
30:675--692, 1981.

\bibitem{Majid92} S.~Majid. \newblock Tannaka-{K}rein {T}heorem for {Q}uasi-{%
H}opf {A}lgebras. \newblock In M.~Gestenhaber and J.~Stasheff, editors, 
\emph{Deformation Theory and Quantum Groups with Applications to
Mathematical Physics, Amherst, {MA}, 1990}, volume 134 of \emph{Contemp.
Math.}, pages 219--232. Amer. Math. Soc., Providence, RI, 1992.

\bibitem{Majid95} S.~Majid. 
\newblock {\em Foundations of {Q}uantum {G}roup
{T}heory}. \newblock Cambridge University Presss, Cambridge, 1995.

\bibitem{Masuoka03} A.~Masuoka. \newblock More {H}omological {A}pproach to {C%
}omposition of {S}ubfactors. \newblock {\em J. Math. Sci., Tokyo},
10(4):599--630, 2003.

\bibitem{Panaite97Stefan} F.~Panaite and D.~\c{S}tefan. \newblock When is
the {C}ategory of {C}omodules a {B}raided {T}ensor {C}ategory? \newblock
\emph{Rev. Roum. Math. Pures Appl.}, 42(1-2):107--119, 1997.

\bibitem{Pareigis77} B.~Pareigis. \newblock Non-{A}dditive {R}ing and {M}%
odule {T}heory {I}. {G}eneral {T}heory of {M}onoids. 
\newblock {\em Publ.
Math. Debrecen}, 24:189--204, 1977.

\bibitem{Pareigis77II} B.~Pareigis. \newblock Non-{A}dditive {R}ing and {M}%
odule {T}heory {II}. {C}-{C}ategories, {C}-{F}unctors and {C}-{M}orphisms. %
\newblock {\em Publ. Math. Debrecen}, 24:351--361, 1977.

\bibitem{Schauenburg98} P.~Schauenburg. \newblock Bialgebras over {N}%
oncommutative {R}ings and a {S}tructure {T}heorem for {H}opf {B}imodules. %
\newblock {\em Appl. Cat. Str.}, 6:193--222, 1998.

\bibitem{Schauenburg03} P.~Schauenburg. \newblock Actions of {M}onoidal {C}%
ategories, and {G}eneralized {H}opf {S}mash {P}roducts. 
\newblock {\em J.
Alg.}, 270:521--563, 2003.

\bibitem{Schauenburg04} P.~Schauenburg. \newblock Hopf-{G}alois and {B}i-{G}%
alois {E}xtensions. \newblock In G.~Janelidze, B.~Pareigis, and W.~Tholen,
editors, \emph{Galois Theory, Hopf Algebras, and Semiabelian Categories},
volume~43 of \emph{Fields Inst. Commun.} AMS, 2004.

\bibitem{Schneider90a} H.-J. Schneider. \newblock Principal {H}omogeneous {S}%
paces for {A}rbitrary {H}opf {A}lgebras. \newblock {\em Israel J. Math.},
72(1-2):167--195, 1990.

\bibitem{Schneider90} H.-J. Schneider. \newblock Representation {T}heory of {%
H}opf {G}alois {E}xtensions. \newblock {\em Israel J. Math.},
72(1-2):196--231, 1990.

\bibitem{Ulbrich81} K.~H. Ulbrich. \newblock Vollgraduierte {A}lgebren. %
\newblock {\em Abh. Math. Sem. Univ. Hamburg}, 51:136--148, 1981.
\end{thebibliography}
\end{document}